\documentclass[11pt]{article}
\newcommand{\documentdate}{7 III 2022}

\usepackage{a4wide,latexsym,amsmath,varioref,graphicx,xcolor}
\title{OFFO minimization algorithms for second-order optimality and their complexity}

\author{
   S. Gratton%
   \thanks{Universit\'e de Toulouse, INP, IRIT, Toulouse, France. Email:
     serge.gratton@enseeiht.fr. Work partially supported by 3IA Artificial and
     Natural Intelligence Toulouse Institute (ANITI), French "Investing for the Future
     - PIA3" program under the Grant agreement ANR-19-PI3A-0004"}
   ~and Ph. L. Toint%
   \thanks{NAXYS, University of Namur, Namur, Belgium. Email:
     philippe.toint@unamur.be. Work partially supported by ANITI.}
}

\newcommand{\beqn}[1]{\begin{equation}\label{#1}}
\newcommand{\eeqn}{\end{equation}}
\newcommand{\req}[1]{(\ref{#1})}
\newcommand{\ms}{\;\;\;\;}
\newcommand{\tim}[1]{\;\; \mbox{#1} \;\;}

\newtheorem{theorem}{Theorem}[section]
\newtheorem{lemma}[theorem]{Lemma}

\newcommand{\numsection}[1]{\section{#1}\setcounter{equation}{0}}
\newtheorem{corollary}{Corallary}
\newcommand{\appnumsection}[1]{\section*{#1}\setcounter{equation}{0}
  \renewcommand{\theequation}{A.\arabic{equation}}
  \renewcommand{\thetheorem}{A.\arabic{theorem}}
  \renewcommand{\thetable}{A.\arabic{table}}
  \renewcommand{\thefigure}{A.\arabic{figure}}
  \renewcommand{\thesection}{A} }
\renewcommand{\theequation}{\arabic{section}.\arabic{equation}}

\newcounter{algo}[section]
\renewcommand{\thealgo}{\thesection.\arabic{algo}}
\newcommand{\llem}[2]{\vspace{\baselineskip} 
\noindent\framebox[\textwidth]{\parbox{0.95\textwidth}{
\begin{lemma} \label{#1} \rm #2 \end{lemma} } } \vspace{\baselineskip} }
\newlength{\thmw}
\newcommand{\algo}[3]{\refstepcounter{algo}
\begin{center}\begin{figure}[htbp]
\framebox[\textwidth]{
\parbox{0.95\textwidth} {\vspace{\topsep}
{\bf Algorithm \thealgo : #2}\label{#1}\\
\vspace*{-\topsep} \mbox{ }\\
{#3} \vspace{\topsep} }}
\end{figure}\end{center}}
\newcommand{\bpr}{{\bf Proof.} \hspace{1.5mm}}
\newcommand{\epr}{\hfill $\Box$ \vspace*{1em}}
\newcommand{\proof}[1]{
\begin{list}{}{
\setlength{\topsep}{0.0pt}
\setlength{\partopsep}{0.0pt}
\setlength{\leftmargin}{0.025\textwidth}
\setlength{\rightmargin}{0.5\leftmargin}
\setlength{\labelwidth}{0.5\leftmargin}
\setlength{\labelsep}{0.25\leftmargin}}
\item \bpr #1 \epr \noindent
\end{list}}
\newcommand{\lthm}[2]{\vspace{\baselineskip} 
\noindent\framebox[\textwidth]{\parbox{0.95\textwidth}{
    \begin{theorem} \label{#1} \rm #2 \end{theorem} } } \vspace{\baselineskip} }
\DeclareMathOperator*{\argmin}{arg\,min}
\newcommand{\pd}[1]{\langle #1 \rangle}
\newcommand{\ii}[1]{\{ 1, \ldots, #1 \}}
\newcommand{\iiz}[1]{\{ 0, \ldots, #1 \}}
\newcommand{\iibe}[2]{\{ #1, \ldots, #2 \}}

\newcommand{\calK}{{\cal K}} 
 
\newcommand{\calO}{{\cal O}} 
\newcommand{\calS}{{\cal S}}

\renewcommand{\Re}{\hbox{I\hskip -2pt R}}
\newcommand{\smallRe}{\hbox{\footnotesize I\hskip -2pt R}}
\newcommand{\bigfrac}[2]{\frac{\displaystyle #1}{\displaystyle #2}}
\newcommand{\bigsum}{\displaystyle \sum}

\newcommand{\sfrac}[2]{{\scriptstyle \frac{#1}{#2}}}
\newcommand{\kap}[1]{\kappa_{\mbox{\tiny #1}}}
\newcommand{\eqdef}{\stackrel{\rm def}{=}}
\newcommand{\mystack}[2]{_{\stackrel{\scriptstyle #1}{\scriptstyle #2}}}
\newcommand{\al}[1]{{\footnotesize{\sf #1}}}

\newcommand{\tal}[1]{{\normalsize {\sf #1}}}
\newcommand{\half}{\sfrac{1}{2}}
\newcommand{\third}{\sfrac{1}{3}}
\newcommand{\quarter}{\sfrac{1}{4}}
\newcommand{\flow}{f_{\rm low}}

\newcommand{\hatw}{\hat{w}}
\newcommand{\hatphi}{\widehat{\phi}}

\DeclareMathOperator*{\average}{average}

\date{\documentdate}

\begin{document}

\maketitle

\begin{abstract}
An Adagrad-inspired class of algorithms for smooth unconstrained optimization
is presented in which the objective function is never evaluated and yet the
gradient norms decrease at least as fast as $\calO(1/\sqrt{k+1})$ while
second-order optimality measures converge to zero at least as fast as
$\calO(1/(k+1)^{1/3})$. This latter rate of convergence is shown to be
essentially sharp and is identical to that known for more standard algorithms
(like trust-region or adaptive-regularization methods) using both function and
derivatives' evaluations.  A related ``divergent stepsize'' method is also
described, whose essentially sharp rate of convergence is slighly inferior. It
is finally discussed how to obtain weaker second-order optimality guarantees
at a (much) reduced computional cost.
\end{abstract}


\numsection{Introduction}

This paper considers an \emph{a priori} unexpected but fundamental and
challenging question: is evaluating the value of the objective function
necessary for obtaining (complexity-wise) efficient minimization algorithms
which find second-order approximate minimizers?  This question arose as a
natural consequence of the somewhat surprising results of
\cite{GratJeraToin22b}, where it was shown that OFFO (i.e.\ Objective-Function
Free Optimization) algorithms\footnote{For which the only source of
  information on the problem at hand is the value of the gradient.}  exist
which converge to first-order points at a global rate which in
order identical to that to well-known methods using both gradient and
objective function evaluations. That these algorithms include the
deterministic version of Adagrad \cite{DuchHazaSing11}, a very popular method
for deep learning applications, was an added bonus and a good motivation.

We show here that, \emph{from the point of view of evaluation complexity
  alone, evaluating the value of the objective function during optimization is
  also unnecessary\footnote{The authors are well aware that this is a
    theoretical statement, as it may be impractical to evaluate derivatives
    without first evaluating the function itself.} for finding approximate
  second-order minimizers} at a (worst-case) cost entirely comparable to that
incurred by familiar and reliable techniques such as second-order trust-region
or adaptive regularization methods.  This conclusion is coherent with that of
\cite{GratJeraToin22b} for first-order points and is obtained by exhibiting an
OFFO algorithm whose global rate of convergence is proved to be
$\calO(1/\sqrt{k+1})$ for the gradients'norm and $\calO(1/(k+1)^{1/3})$ for
second-order measures.  The new \al{ASTR2} algorithm is of the
adaptively scaled trust-region type, as those studied in
\cite{GratJeraToin22b}.  The key difference is that it now hinges on a scaling
technique which depends on second-order information, when relevant.

The paper is organized as follows. Section~\ref{algo-s} presents the new
\al{ASTR2} class of algorithms and discusses some of its scaling-independent
properties. The complexity analysis of a first,
Adagrad-like, subclass of \al{ASTR2} is then presented in
Section~\ref{adag-s}. Another subclass of interest is also considered and
analyzed in Section~\ref{ds-s}.  Section~\ref{weak2nd-s} discusses how weaker
optimality conditions may be guaranteed by the \al{ASTR2} algorithms at
significantly reduced computational cost. Conclusions and perspectives are finally
presented in Section~\ref{concl-s}

\numsection{The \tal{ASTR2} class of minimization methods}\label{algo-s}

\subsection{Approximate first- and second-order optimality}

We consider the nonlinear unconstrained optimization problem
\beqn{problem}
\min_{x\in\smallRe^n} f(x)
\eeqn
where $f$ is a function from $\Re^n$ to $\Re$. More precisely, we assume
that

\vspace*{2mm}
\noindent
\begin{description}
\item[AS.1:] the objective function $f(x)$ is twice continuously differentiable;
\item[AS.2:] 
  its gradient $g(x)\eqdef \nabla_x^1f(x)$ and Hessian $H(x) \eqdef \nabla_x^2f(x)$ are Lipschitz continuous with
   Lipschitz constant $L_1$ and $L_2$, respectively, that is
   \[
   \|g(x)-g(y)\| \leq L_1 \|x-y\|
   \tim{ and }
   \|H(x)-H(y)\| \leq L_2 \|x-y\|
   \]
   for all $x,y\in \Re^n$;
\item[AS.3:] there exists a constant $\flow$ such that $f(x)\ge \flow$ for all $x\in\Re^n$.
\end{description}

As our purpose is to find approximate first- and second-order minimizers, we
need to clarify these concepts. In this paper we choose to follow the ``
strong $\phi$'' concept of optimality discussed in \cite{CartGoulToin20c,CartGoulToin22b} or
\cite[Chapters~12--14]{CartGoulToin22}.  It is based on the quantity
\beqn{phi-def}
\phi_{f,2}^\delta(x) = f(x) - \min_{\|d\|\le \delta} T_{f,2}(x,d),
\eeqn
where $T_{f,2}(x,d)$ is the second-order Taylor expansion of $f$ at $x$, that
is
\[
T_{f,1}(x,d) = f(x) + g(x)^Td
\tim{ and }
T_{f,2}(x,d) = f(x) + g(x)^Td + \half d^TH(x)d.
\]
Observe that $\phi_{f,j}^\delta(x)$ is interpreted as the maximum decrease of
the local $j$-th order Taylor model of the objective function $f$ at $x$, within a
ball of radius $\delta$.  Importantly for our present purposes, the
evaluation of $\phi_{f,2}^\delta(x)$ \emph{does not require the
evaluation of $f(x)$,} as it can be rewritten as
\beqn{phik-2}
\phi_{f,2}^\delta(x) = \max_{\|d\|\le \delta} -\Big( g(x)^Td + \half d^TH(x)d \Big).
\eeqn
The next result recalls the link between
the $\phi$ optimality measure and the more standard ones.

\llem{exact-optim}{\cite[Theorems 12.1.4 and 12.1.6]{CartGoulToin22}
Suppose that $f$ is twice continuously differentiable. Then
\begin{itemize}
\item[(i)] for any $\delta>0$ and any $x\in \Re^n$,  we have that
\beqn{unco-1st-order}
\|g(x)\|=\frac{\phi_{f,1}^{\delta}(x)}{\delta},
\eeqn
and so $\phi_{f,1}^{\delta}(x) = 0$ if and only if $g(x)=0$;
\item[(ii)] we have that
\[
\phi_{f,2}^{\delta}(x) = 0 \tim{for some $\delta>0$, then }
 g(x)=0\,\tim{and}\,\lambda_{\min}[H(x)] \geq 0,
\]
and so any such $x$ is a first- and second-order minimizer;
\item[(iii)] if $\phi_{f,1}^{\delta_1}(x) \leq \epsilon_1 \, \delta_1$ 
(and so \req{optim} holds with $j=1$), then
$\|g(x)\|\leq \epsilon_1$; 
\item[(iv)] if $\phi_{f,2}^\delta(x) \leq \half \epsilon_2\delta^2$, 
then 
$\lambda_{\min}[H(x)] \geq -\epsilon_2$ 
(and so \req{optim} holds for $j=2$) and
$\|g(x)\|\leq \delta\kappa(x)\sqrt{\epsilon_2}$, 
where $\kappa(x)$ depends on (the eigenvalues of)
$H(x)$.
\end{itemize}
}
  
\noindent
Note also that computing $\phi_{f,1}^\delta(x)$  simply results from
\req{unco-1st-order} and that, in particular, $\phi_{f,1}^1(x)=\|g(x)\|$.
Computing $\phi_{f,2}^\delta(x)$ is a standard Euclidean
trust-region step calculation (see \cite[Chapter~7]{ConnGoulToin00}, for
instance).

For $j\in\{1,2\}$, we then say that an iterate $x_k$ is an 
$\epsilon$-approximate minimizer if
\beqn{optim}
\phi_{f,j}^\delta(x_k) \le  \epsilon_j \frac{\delta^j}{j}
\ms
\tim{for some } \delta \in (0,1] \tim{ and all } 1\leq i \leq j,
  \eeqn
where $\epsilon = (\epsilon_1, \ldots, \epsilon_j)$. There are two ways to
express how fast an algorithm tends to such points in the worst case.  The
first (the ``$\epsilon$-orders'') is to assume
$\epsilon$ is given and then give a bound on the maximum number of iterations and
evaluations that are needed to satisfy \req{optim}.  In this paper we focus on
the second (the ``$k$-orders''), where one instead gives an upper
bound\footnote{Converging to zero.} on
$\phi_{f,j}^\delta(x_k)$ as a function of $k$ (for specified $j$ and $\delta$).

\subsection{The \tal{ASTR2} class}

After these preliminaries, we now introduce the new \al{ASTR2} class of
algorithms.  Methods in this class are of ``adaptively scaled trust-region''
type, a term we now briefly explain.  Classical trust-region algorithms (see
\cite{ConnGoulToin00} for an in-depth coverage or \cite{Yuan15} for a more
recent survey) are iterative. At each iteration, they define a local model of
the objective function which is deemed trustable within the ``trust region'',
a ball of given radius centered at the current iterate. A step and
corresponding trial point are then computed by (possibly approximately)
minimizing this model in the trust region.  The objective function value is
then computed at the trial point, and this point is accepted as the new
iterate if the ratio of the achieved reduction in the objective function to
that predicted by the model is sufficiently large.  The radius of the trust
region is then updated using the value of this ratio.  As is clear from this
description, these methods are intrinsically dependent of the evaluation of
the objective function, and therefore not suited to our Objective-Function
Free Optimization (OFFO) context. Here we follow \cite{GratJeraToin22b} in
interpreting the mechanism designed for the Adagrad methods
\cite{DuchHazaSing11} as an alternative trust-region design not using function
evaluations.  In this interpretation, the trial point is always accepted and
the trust-region radius is determined by the gradient sizes, in a manner
reminiscent also of \cite{FanYuan01}. In this approach, one uses {\em scaling
  factors} to determine the radius (hence the name of Adaptively Scaled Trust
Region).  Given these factors, we may then state the \al{ASTR2} class of
algorithms as shown \vpageref{ASTR2}. This algorithm involves requirements on
the step which are standard (and practical) for trust-region methods.
  
\algo{ASTR2}{\tal{ASTR2}}
{
\begin{description}
\item[Step 0: Initialization. ]
A starting point $x_0$ is given. The constants $\tau,\chi \in (0,1]$ and
  $\xi \geq 1$ are also given.
  Set $k=0$.
\item[Step 1: Compute derivatives. ]
  Compute $g_k = g(x_k)$ and $H_k = H(x_k)$, as well as
  $\phi_k \eqdef \phi_{f,2}^1(x_k)$ and $\hatphi_k \eqdef \min[\phi_k,\xi]$.
%  If $\|g_k\|\leq \epsilon_1$ and $\hatphi_k\leq \half \epsilon_2$, terminate.
\item[Step 2: Define the trust-region radii. ]
  Set
  \beqn{Delta-def}
  \Delta^L_k  = \frac{\|g_k\|}{w^L_k}
  \tim{ and }
  \Delta^Q_k  = \frac{\hatphi_k}{w^Q_k}
  \eeqn
  where $w^L_k=w^L(x_0,\ldots,x_k)$ and $w^Q_k=w^Q(x_0,\ldots,x_k)$.
\item[Step 3: Step computation.]
  If
  \beqn{LQ-cond}
  \|g_k\|^2 \geq \hatphi_k^3,
  \eeqn
  then set 
  \beqn{sL-def}
  s_k = s^L_k = -\frac{g_k}{w^L_k}.
  \eeqn
  Otherwise, set $s_k= s^Q_k$, where $s^Q_k$ is such that
  \beqn{sQ-def}
  \|s^Q_k\| \leq \Delta^Q_k
  \tim{ and }
  \Delta q_k \ge \tau \max\big[\Delta q_k^C,\Delta q_k^E\Big]
  \eeqn
  where
  \beqn{C-decr}
  \Delta q_k^C
  = \max_{\mystack{\alpha\geq0}{\alpha\|g_k\|\leq\Delta_k^Q}}
        \Big[f(x_k)-T_{f,2}(x_k,-\alpha g_k )\Big]
  \eeqn
  and
  \beqn{E-decr}
  \Delta q_k^E 
  = \max_{\mystack{\alpha\geq0}{\alpha\leq\Delta_k^Q}}
        \Big[f(x_k)-T_{f,2}(x_k, \alpha u_k )\Big]
  \eeqn
  with $u_k$ satisfying
  \beqn{uk-conds}
  u_k^TH_ku_k \leq \chi \, \lambda_{\min}[H_k], \ms u_k^Tg_k \leq 0
  \tim{ and } \|u_k\|=1.
  \eeqn
\item[Step 4: New iterate.]
  Define
  \beqn{xupdate-a}
  x_{k+1} = x_k + s_k,
  \eeqn
  increment $k$ by one and return to Step~1.
\end{description}
}

\noindent
A few additional  comments on this algorithm are now in order.
\begin{enumerate}
\item The algorithms in the \al{ASTR2} class belong to the OFFO framework: the objective
  function is never evaluated (remember that $\phi_{f,j}^1(x)$ can be computed
  without any such evaluation, the same being obviously true for $\Delta q_k$,
  $\Delta q_k^C$ and $\Delta q_k^E$).
\item Given our focus on $k$-orders of convergence, the algorithm does not
  include a termination criterion.  It is however easy, should one be interested
  in $\epsilon$-orders instead, to test \req{optim} for $\delta=1$ and the considered
  $\epsilon_1$ and $\epsilon_2$ at the end of Step~1, and then terminate if this condition holds.
\item Despite their somewhat daunting statements, conditions
  \req{sQ-def}--\req{uk-conds} are relatively mild and have been extensively
  used for standard trust-region algorithms, both in theory and practice.
  Condition \req{C-decr} defines the so-called ``Cauchy decrease'', which is
  the decrease achievable on the quadratic model $T_{f,2}(x_k,s)$ in the
  steepest descent direction \cite[Section~6.3.2]{ConnGoulToin00}. Conditions
  \req{E-decr} and \req{uk-conds} define the ``eigen-point decrease'', which
  is that achievable along $u_k$, a ($\chi$-approximate) eigenvector
  associated with the smallest Hessian eigenvalue \cite[Section~6.6]{ConnGoulToin00}.
  We discuss in Section~\ref{weak2nd-s} how they can be ensured in practice,
  possibly approximately, for instance by the GLTR algorithm \cite{GoulLuciRomaToin99}.
\item The computation of $\phi_k$ can be reused to compute $s^Q_k$, should it
  be necessary. If $\Delta_k>1$, the model minimization may be pursued beyond
  the boundary of the unit ball. If $\Delta_k<1$, backtracking is also possible
  \cite[Section~10.3.2]{ConnGoulToin00}.
\item Note that \emph{two} scaling factors are updated from iteration to iteration:
  one for first-order models and one for second-order ones.  It does indeed
  make sense to trust these two types of models in region of different sizes,
  as Taylor's theory suggests second-order models may be reasonably accurate
  in larger neighbourhoods.
\item A ``componentwise'' version where the trust region is defined in the
  $\|\cdot\|_\infty$ norm is possible with
  \[
  \phi_{i,k} = \max\left[\phi_k, -
     \min_{|\alpha| \le 1}\left(\alpha g_{i,k} +\half\alpha^2
      [H_k]_{i,i}\right)\right]
  \]
  and
  \[
  \Delta^L_{i,k}  = \frac{|g_{i,k}|}{w^L_{i,k}}
  \tim{ and }
  \Delta^Q_{i,k}  = \frac{\min[\xi, \phi_{i,k}]}{w^Q_{i,k}}.
  \]
  We will not explicitly consider this variant to keep our notations
  reasonably simple.
\end{enumerate}

Our assumption that the gradient and Hessian are Lipschitz continuous (AS.2) ensures the
following standard result.

\llem{lemma:Lipschitz}{\cite{BirgGardMartSantToin17} or \cite[Theorem~A.8.3]{CartGoulToin22}
Suppose that AS.1 and AS.2 hold.  Then
\beqn{lin-Lip}
f(x_k+s^L_k)- f(x_k) \leq \pd{g_k,s^L_k} + \frac{L_1}{2} \|s^L_k\|^2
\eeqn
and
\beqn{quad-Lip}
f(x_k+s^Q_k) - f(x_k) \leq -\Delta q_k  + \frac{L_2}{6}\|s^Q_k\|^3.
\eeqn
}

\noindent
The first step in analyzing the convergence of the \al{ASTR2} algorithm is to
derive bounds on the objective function's change from iteration to iteration,
depending on which step (linear with $s_k=s_k^L$, or quadratic with $s_k=s_k^Q$) is chosen. We start by a few
auxiliary results on the relations between first- and second-order optimality
measures.

\llem{phi-convex}{Suppose that $H$ is an $n\times n$ symmetric positive semi-definite
  matrix and $g \in \Re^n$, and consider the (convex) quadratic
  $q(d) = \pd{g,d} + \half \pd{d,Hd}$.  Then
  \beqn{phi-ng-convex}
  \phi_{q,2}^1(0) = \left|\min_{\|d\|\leq 1}q(d)\right| \leq \|g\|.
  \eeqn
}

\proof{
  From the definition of the gradient, we have that
  \[
  \|g\| = \left|\min_{\|d\|\leq 1}\pd{g,d}\right|.
  \]
  But $\pd{g,d}$ defines the supporting hyperplane of $q(d)$ at $d=0$
  and thus the convexity of $q$ implies that $q(d) \geq \pd{g,d}$ for all $d$.
  Hence
  \[
  \left|\min_{\|d\|\leq 1}q(d)\right|
  \leq\left|\min_{\|d\|\leq 1}\pd{g,d}\right|
  \]
  and \req{phi-ng-convex} follows.
}  %epr

\llem{phi-lambga-g}{Suppose that
  \beqn{near-convex}
  0 < \eta_k\leq \half \phi_k
  \eeqn
  where
  \beqn{eta-def}
  \eta_k \eqdef \min\bigg(0,-\lambda_{\min}[H_k]\bigg).
  \eeqn
  Then
  \beqn{phi-ng}
  \half \phi_k \leq \|g_k\|.
  \eeqn
}

\proof{
  Observe first that \req{near-convex} implies that $\lambda_{\min}[H_k] < 0$
  and $\eta_k = \left|\lambda_{\min}[H_k]\right|$.
  Let $d_k$ be a solution of the optimization problem defining $\phi_k$, i.e.,
  \[
  d_k = \argmin_{\|d\|\leq 1}T_{f,2}(x_k,d),
  \]
  so that $\phi_k = f_k-T_{f,2}(x_k,d_k)$.
  Since $\lambda_{\min}[H_k] < 0$, it is known from trust-region theory
  \cite[Corollary.2.2]{ConnGoulToin00} that $d_k$ may be chosen such that
  $\|d_k\| = 1$. Now define
  \[
  q_0(d) \eqdef \pd{g_k,d} + \half \pd{d,(H_k-\lambda_{\min}[H_k] I)d} 
    = T_{f,2}(x_k,d) - f_k + \eta_k \|d\|^2 
  \]
  and note that $q_0(d)$ is convex by construction. Then, at $d_k$,
  \[
  q_0(d_k) = -\phi_k + \eta_k
  \]
  and \req{near-convex} implies that $q_0(d_k) < 0$.  Moreover,
  \[
  \half \phi_k \leq - q_0(d_k) \leq -\pd{g_k,d_k} \leq \|g_k\|,
  \]
  where we used the convexity of $q_0$ to deduce the the first inequality, and
  Cauchy-Schwarz with $\|d_k\| \leq 1$ to derive the second. This proves \req{phi-ng}.
} % epr

\noindent
Using these results, we may now prove a crucial property on objective function
change. For this purpose, we partition the iterations in two sets, depending
which type of step is chosen, that is
\[
\calK^L = \{ k \geq 0 \mid s_k = s_k^L \}
\tim{ and }
\calK^Q = \{ k \geq 0 \mid s_k = s_k^Q \}.
\]

\llem{lemma:fdecr}{
Suppose that AS.1 and AS.2 hold. Then
  \beqn{ffdecr0}
  f_{k+1}-f_k
  \leq - \frac{\|g_k\|^2}{w^L_k} + \frac{L_1}{2}\frac{\|g_k\|^2}{(w^L_k)^2}
  \tim{ for } k \in \calK^L
  \eeqn
  and
  \beqn{ffdecr1}
  f_{k+1}-f_k \leq - \frac{\tau}{4\xi}
  \min\left[\frac{1}{2(1+L_1)},\frac{1}{w^Q_k},\frac{1}{(w^Q_k)^2}\right]\,\hatphi_k^3
       + \frac{L_2}{6} \frac{\hatphi_k^3}{(w^Q_k)^3}
  \tim{ for } k \in \calK^Q.
  \eeqn
}

\proof{
  Suppose first that $s_k=s^L_k$. Then \req{lin-Lip}, \req{sL-def}
  and \req{Delta-def} ensure that
  \beqn{decr1}
  f_{k+1}-f_k
  \leq -\frac{\|g_k\|^2}{w^L_k} + \frac{L_1}{2}(\Delta^L_k)^2
  = -\frac{\|g_k\|^2}{w^L_k} + \frac{L_1}{2}\frac{\|g_k\|^2}{(w^L_k)^2},
  \eeqn
  giving \req{ffdecr0}.
  
  Suppose now that $s_k= s^Q_k$, i.e.\ $k\in\calK^Q$. Then, because of
  \req{sQ-def}--\req{uk-conds}, the
  decrease $\Delta q_k$ in the quadratic model $T_{f,2}(x_k,s)$ at $s_k$ is at least a
  fraction $\tau$ of the maximum of the Cauchy and eigen-point decreases given
  by \req{C-decr} and \req{E-decr}. Standard trust-region theory
  (see \cite[Lemmas~6.3.2 and 6.6.1]{ConnGoulToin00} for instance)
  then ensures that, for possibly non-convex $T_{f,2}(x_k,s)$,
  \[
  \begin{array}{lcl}
  \Delta q_k
  & \geq & \tau\max\left[\bigfrac{1}{2}\min\left(\bigfrac{\|g_k\|^2}{1+\|H_k\|},\|g_k\|\Delta^Q_k\right),
                \bigfrac{\eta_k}{2} (\Delta^Q_k)^2\right] \\*[3ex]
  & \geq & \bigfrac{\tau}{2}\max\left[\min\left(\bigfrac{\|g_k\|^2}{1+L_1},\bigfrac{\|g_k\|\hatphi_k}{w^Q_k}\right),
                \bigfrac{\eta_k\hatphi_k^2}{(w^Q_k)^2}\right]
  \end{array}
  \]
  where we used the bound $\|H_k\| \leq L_1$ and \req{Delta-def} to derive the
  last inequality.
  If $\eta_k\leq \half \phi_k$, then, using Lemma~\ref{phi-lambga-g} and the
  inequality $\phi_k\geq \hatphi_k$,
  \[
  \Delta q_k
  \geq \bigfrac{\tau}{2}\min\left(\bigfrac{\|g_k\|^2}{1+L_1},\bigfrac{\|g_k\|\hatphi_k}{w^Q_k}\right)
  \geq \bigfrac{\tau}{2}\min\left(\bigfrac{(\half\hatphi_k)^2}{1+L_1},\bigfrac{(\half\hatphi_k)\hatphi_k}{w^Q_k}\right).
  \]
  Now $\hatphi_k^3 \leq \xi \hatphi_k^2$ and thus
  \beqn{TR1}
  \Delta q_k
  \geq
  \bigfrac{\tau}{2}\min\left(\bigfrac{\hatphi_k^3}{4\xi(1+L_1)},
  \bigfrac{\hatphi_k^3}{2\xi w^Q_k}\right).
  \eeqn
  If instead $\eta_k > \half \phi_k\geq \half \hatphi_k$, then
  \beqn{TR2}
  \Delta q_k
  \geq \bigfrac{\tau}{2}\bigfrac{\eta_k\hatphi_k^2}{(w^Q_k)^2}
  \geq \bigfrac{\tau}{2}\bigfrac{(\half\hatphi_k)\hatphi_k^2}{(w^Q_k)^2}.
  \eeqn
  Given that, if $k\in\calK^Q$, $\|s_k\|\leq \Delta_k^Q=
  \hatphi_k/w_k^Q$, we deduce \req{ffdecr1} from \req{quad-Lip}, \req{TR1} and \req{TR2}.
}  % epr

\noindent
Observe that neither \req{ffdecr0} nor \req{ffdecr1} guarantees that the
objective function values are monotonically decreasing.

\numsection{An Adagrad-like algorithm for second-order optimality}\label{adag-s}

We first consider a choice of scaling factors directly inspired by the
Adagrad algorithm \cite{DuchHazaSing11} and assume that, for some
$\varsigma > 0$, $\mu,\nu \in (0,1)$, $\vartheta_L,\vartheta_Q\in(0,1]$ and all $k\geq0$,
\beqn{w-adag-L}
w_k^L \in [\vartheta_L \hatw_k^L, \hatw_k^L]
\tim{ where }
\hatw^L_k = \left(\varsigma + \sum_{\mystack{\ell=0}{\ell\in\calK^L}}^k
\|g_k\|^2\right)^\mu
\eeqn
and
\beqn{w-adag-Q}
w_k^Q \in [\vartheta_Q \hatw_k^Q,\hatw_k^Q]
\tim{ where }
\hatw^Q_k = \left(\varsigma + \sum_{\mystack{\ell=0}{\ell\in\calK^Q}}^k \hatphi_k^3\right)^\nu.
\eeqn
Note that selecting the parameters $\vartheta_L$ and $\vartheta_Q$ strictly
less than one allows the scaling factors $w_k^L$ and $w_k^Q$ to be chosen in
an interval at each iteration without any monotonicity.

We now present a two technical lemmas which will be
necessary in our analysis.  The first states useful results for a specific
class of inequalities.

\llem{techsum}{
Let $a\geq \half\varsigma$ and $b\geq \half\varsigma$. Suppose
that, for some $\theta_a\geq 1$,
$\theta_b\geq 1$, $\theta\geq 0$, $\mu\in(0,1)$,  and $\nu\in(0,\third)$
\beqn{techsum:hyp}
a^{1-\mu} + b^{1-2\nu} \leq \theta_a A(a) + \theta_b B(b) + \theta 
\eeqn
where $A(a)$ and $B(b)$ are given, as a function of $\mu$ and $\nu$, by\\*[1ex]
\centerline{
%     \begin{center}
     \begin{tabular}{|l|c|c|c|}
     \hline
           & $\mu<\half$  & $\mu = \half$ & $\mu>\half$ \\
     \hline
       $A(a)$ & $a^{1-2\mu}$ & $\log(2a)$ & 0 \\
     \hline
     \end{tabular}
     ~and~
     \begin{tabular}{|l|c|c|c|}
     \hline
              & $\nu<\third$ & $\nu = \third$   & $\nu>\third$ \\
     \hline
       $B(k)$ & $b^{1-3\nu}$ &   $\log(2b)$   &   0   \\
      \hline
     \end{tabular} .
     %     \end{center}
     }\\*[1ex]
     Then there exists positive constants $\kappa_a$ and $\kappa_b$ only
     depending on $\theta_a$, $\theta_b$, $\theta$, $\mu$ and $\nu$ such that
     \beqn{akbk-bound}
     a\leq \kappa_a  \tim{ and } b \leq \kappa_b.
     \eeqn
}

\proof{
This result is proved by comparing the value of the left- and right-hand sides
for possibly large $a$ and $b$.  The details are given in Lemmas~\ref{techa}--\ref{techg}
in appendix, whose results are then combined as shown in Table~\ref{cases}.
  \begin{table}[htb]
    \begin{center}
      \begin{tabular}{|l|c|c|c|}
        \hline
                    & $\mu<\half$  & $\mu = \half$ & $\mu>\half$ \\
        \hline
        $\nu<\third$&   Lemma~\ref{techc}  &   Lemma~\ref{techf}   &   Lemma~\ref{techd}  \\
        $\nu=\third$&   Lemma~\ref{techf}  &   Lemma~\ref{techg}   &   Lemma~\ref{techf}  \\
        $\nu>\third$&   Lemma~\ref{techd}  &   Lemma~\ref{techf}   &   Lemma~\ref{techa}  \\
        \hline
      \end{tabular}
      \caption{\label{cases}Lemmas for combinations of $\mu$ and $\nu$}
    \end{center}
  \end{table}
The details of the constants $\kappa_a$ and $\kappa_b$ for the various cases
are explicitly given in the statements of the relevant lemmas.
}

\noindent
The second auxiliary result is a bound extracted from
\cite{GratJeraToin22b} (see also \cite{DefoBottBachUsun20,WardWuBott19} for
the case $\alpha=1$).

\llem{gen:series}{Let $\{c_k\}$ be a non-negative sequence,
  $\varsigma>0$, $\alpha > 0$, $\nu \geq 0$ and define, for each $k \geq 0$,
	$d_k = \sum_{j=0}^k c_j$.  If $\alpha \neq 1 $, then
\beqn{allalpha series-bound}
\sum_{j=0}^k \frac{c_j}{(\varsigma+d_j)^{\alpha}}
\le \frac{1}{(1-\alpha)} ( (\varsigma + d_k)^{1-  \alpha} - \varsigma^{1-  \alpha} ).
\eeqn
Otherwise,
\beqn{alphasup1series-bound}
\sum_{j=0}^k\frac{c_j}{(\varsigma+d_j)}
\le  \log\left(\frac{\varsigma + d_k}{\varsigma} \right).
\eeqn
}

\noindent
Note that, if $\alpha > 1$, then the bound \req{allalpha series-bound} can be
rewritten as
\[
\sum_{j=0}^k  \frac{ c_j}{(\varsigma+d_j)^{\alpha}}
\le \frac{1}{\alpha-1} \Big( \varsigma^{1-\alpha} -(\varsigma + d_k)^{1-\alpha}\Big),
\]
whose right-hand side is positive.

\noindent
Armed with the above results, we are now in position to specify particular
choices of the scaling factors $w_k$ and derive the convergence properties
of the resulting variants of \al{ASTR2}.

\lthm{theorem:better-complexity}{Suppose that AS.1--AS.3 hold and that the
  \al{ASTR2} algorithm is applied to problem \req{problem}, where $w^L_k$ 
  and $w^Q_k$ are given by \req{w-adag-L} and \req{w-adag-Q}, respectively. Then there exists a positive constant
  $\kap{ASTR2}$ only depending on the problem-related
  quantities $x_0$, $\flow$, $L_1$ and $L_2$  and on the algorithmic parameters
  $\varsigma$, $\tau$, $\xi$, $\mu$ and $\nu$ such that
  \beqn{avg-orders}
  \average_{j\in\iiz{k}} \|g_j\|^2   \leq \frac{\kap{ASTR2}}{k+1}
  \tim{ and }
  \average_{j\in\iiz{k}} \hatphi_j^3 \leq \frac{\kap{ASTR2}}{k+1},
  \eeqn
  and therefore that
  \beqn{min-orders}
  \min_{j\in\iiz{k}} \|g_j\|   \leq \frac{\kap{ASTR2}}{(k+1)^\half}
  \tim{ and }
  \min_{j\in\iiz{k}} \hatphi_j \leq \frac{\kap{ASTR2}}{(k+1)^\third}.
  \eeqn
}

\proof{To simplify notations in the proof,  define
  \beqn{akbk-def}
  a_k = 2\sum_{\mystack{j=0}{j\in\calK^L}}^k \|g_j\|^2
  \tim{ and }
  b_k =  2\sum_{\mystack{j=0}{j\in\calK^Q}}^k \hatphi_k^3.
  \eeqn
  Consider first an iteration index $j\in\calK^L$
  Then \req{ffdecr0} (expressed for for $j\ge0$), \req{w-adag-L} and the inequality
  $\tau\leq 1$ give that
  \beqn{decrj1-n}
  f(x_{j+1})-f(x_j)
  \le - \bigfrac{\tau}{2} \bigfrac{\|g_j\|^2}{w^L_j}
            + \bigfrac{L_1}{2\vartheta_L^2} \bigfrac{\|g_j\|^2}{(w^L_j)^2}
  \le - \bigfrac{\tau}{2} \bigfrac{\|g_j\|^2}{(\varsigma+\half a_j)^\mu}
            + \bigfrac{L_1}{2\vartheta_L^2} \bigfrac{\|g_j\|^2}{(\varsigma+\half a_j)^{2\mu}}.
  \eeqn
  Suppose now that $j\in\calK^Q$.  Then \req{ffdecr1} and \req{w-adag-Q} imply that
  \beqn{decrj2-n}
  f_{j+1}-f_j
  \leq -
  \frac{\tau}{4\xi}\min\left[\frac{\hatphi_j^3}{2(1+L_1)},
    \frac{\hatphi_j^3}{(\varsigma+\half b_j)^\nu}
    \frac{\hatphi_j^3}{(\varsigma+\half b_j)^{2\nu}}
    \right]
       + \frac{L_2}{6\vartheta_Q^3} \frac{\hatphi_j^3}{(\varsigma+\half b_j)^{3\nu}},
  \eeqn
  Suppose now that
  \beqn{large1}
  a_j > 2\varsigma
  \tim{ and}
  b_j > \max\left[1, 2\varsigma,\Big(2(1+L_1)\Big)^\sfrac{1}{\nu}\right],
  \eeqn
  which implies that
  \[
  w_j^L \leq a_j^\mu,
  \ms
  w_j^Q \leq b_j^\nu
  \tim{ and }
  2(1+L_1) \leq  b_j^\nu.
  \]
  Then combining \req{decrj1-n} and \req{decrj2-n}, the inequality
  $\xi\geq 1$ and AS.3, we deduce that, for all $k\geq 0$,
  \[
  f(x_0)-\flow
  \geq \frac{\tau}{4\xi} \left[
    \bigsum_{\mystack{j=0}{j\in\calK^L}}^k  \bigfrac{\|g_j\|^2}{a_j^\mu} +
    \bigsum_{\mystack{j=0}{j\in\calK^q}}^k \bigfrac{\hatphi_j^3}{b_j^{2\nu}}\right]
    - \bigfrac{L_1}{2\vartheta_L^2} \bigsum_{\mystack{j=0}{j\in\calK^L}}^k \bigfrac{\|g_j\|^2}{(w_k^L)^2}
    - \bigfrac{L_2}{6\vartheta_Q^3} \bigsum_{\mystack{j=0}{j\in\calK^Q}}^k \bigfrac{\hatphi_j^3}{(w_k^Q)^3}.
  \]
  But, by definition, $a_j\leq a_k$ and $b_j\leq b_k$ for $j\leq k$, and thus,
  for all $k\geq 0$,
  \beqn{glob-bnd}
  a_k^{1-\mu} + b_k^{1-2\nu}
  \leq \bigfrac{4\xi(f(x_0)-\flow)}{\tau}
    + \bigfrac{2\xi L_1}{ \tau\vartheta_L^2} \bigsum_{\mystack{j=0}{j\in\calK^L}}^k \bigfrac{\|g_j\|^2}{(w_k^L)^2}
    + \bigfrac{2\xi L_2}{3\tau\vartheta_Q^3} \bigsum_{\mystack{j=0}{j\in\calK^Q}}^k \bigfrac{\hatphi_j^3}{(w_k^Q)^3}.
  \eeqn
  We now have to bound the last two terms on the right-hand side of
  \req{glob-bnd}. Using \req{w-adag-L} and
  Lemma~\ref{gen:series} with $\{c_k\} =   \{\|g_k\|^2\}_{k\in\calK^L}$  and
  $\alpha = 2\mu$, gives that 
  \beqn{nolog21-n}
  \bigsum_{\mystack{j=0}{k\in\calK^L}}^k \bigfrac{\|g_j\|^2}{(w^L_k)^2}
  \le \frac{1}{\vartheta_L^2(1-2\mu)} \left( \Big(\varsigma + \sum_{\mystack{j=0}{k\in\calK^L}}^k\|g_k\|^2\Big)^{1-2\mu} - \varsigma^{1-2\mu} \right)
  \le \frac{a_k^{1-2\mu}}{\vartheta_L^2(1-2\mu)}
  \eeqn
   if $\mu < \half$, and
  \beqn{log2-n}
  \bigsum_{\mystack{j=0}{k\in\calK^L}}^k \bigfrac{\|g_j\|^2}{(w^L_j)^2}
  \le \frac{1}{\vartheta_L^2}\log\left(\frac{\varsigma+ \sum_{j=0,k\in\calK^L}^k\|g_k\|^2}{\varsigma}\right)
  \le \frac{1}{\vartheta_L^2}\log\left(\frac{\varsigma+a_k}{\varsigma}\right)
  \eeqn
  if $\mu = \half$ and
  \beqn{nolog22-n}
  \bigsum_{\mystack{j=0}{k\in\calK^L}}^k \bigfrac{\|g_j\|^2}{(w^L_k)^2}
  \le \frac{1}{\vartheta_L^2(2\mu-1)}
  \left(\varsigma^{1-2\mu} - \big(\varsigma + \sum_{\mystack{j=0}{k\in\calK^L}}^k\|g_k\|^2\big)^{1-2\mu}\right)
  \le \frac{\varsigma^{1-2\mu}}{\vartheta_L^2(2\mu-1)}
  \eeqn
  if $\mu > \half$.  Similarly, using \req{w-adag-Q} and Lemma~\ref{gen:series} with $\{c_k\} =
  \{\hatphi_k^3\}_{k\in\calK^Q}$ and  $\alpha = 3\nu$ yields that
  \beqn{nolog31-n}
  \bigsum_{\mystack{j=0}{k\in\calK^Q}}^k \bigfrac{\phi_j^3}{(w^Q_k)^3}
  \le \frac{1}{\vartheta_Q^3(1-3\nu)} \left( \big(\varsigma +
  \sum_{\mystack{j=0}{k\in\calK^Q}}^k\hatphi_k^3\big)^{1-3\nu} - \varsigma^{1-3\nu} \right)
  \le \frac{b_k^{1-3\nu}}{\vartheta_Q^3(1-3\nu)}
  \eeqn
  if $\nu < \third$, 
  \beqn{log3-n}
  \bigsum_{\mystack{j=0}{k\in\calK^Q}}^k \bigfrac{\hatphi_j^3}{(w^Q_j)^3}
  \le \frac{1}{\vartheta_Q^3}\log\left(\frac{\varsigma+ \sum_{j=0,k\in\calK^Q}^k\hatphi_k^3}{\varsigma}\right)
  \le \frac{1}{\vartheta_Q^3}\log\left(\frac{\varsigma+b_k}{\varsigma}\right)
  \eeqn
  if $\nu = \third$, and
  \beqn{nolog32-n}
  \bigsum_{\mystack{j=0}{k\in\calK^Q}}^k \bigfrac{\hatphi_j^3}{(w^Q_k)^3}
  \le \frac{1}{\vartheta_Q^3(3\nu-1)} \left( \varsigma^{1-3\nu}
  - \big(\varsigma + \sum_{\mystack{j=0}{k\in\calK^Q}}^k\hatphi_j^3\big)\right)
  \le \frac{\varsigma^{1-3\nu}}{\vartheta_Q^3(3\nu-1)}
  \eeqn
  if $\nu > \third$.
  Moreover, unless $a_k< 1$, the argument of the logarithm in the
  right-hand side of \req{log2-n} satisfies
  \beqn{arglog2}
  1 \leq \frac{\varsigma+a_k}{\varsigma} \leq 1+a_k \leq 2a_k.
  \eeqn
  Similarly, unless $b_k< 1$, the argument of the logarithm in the
  right-hand side of \req{log3-n} satisfies
  \beqn{arglog3}
  1 \leq \frac{\varsigma+b_k}{\varsigma} \leq 1+b_k \leq 2b_k.
  \eeqn
  Moreover, we may assume, without loss of generality, that $L_1$ and $L_2$
  are large enough to ensure that
  \[
  2\xi L_1 \geq \tau\vartheta_L^2
  \tim{ and }
  2\xi L_2 \geq 3 \tau\vartheta_Q^3.
  \]
  Because of these observations and since \req{glob-bnd} together with one
  of \req{nolog21-n}--\req{nolog22-n} and one of
  \req{nolog31-n}-\req{nolog32-n} has the form of condition 
  \req{techsum:hyp}, we may then apply Lemma~\ref{techsum} for each $k\geq0$
  with $a=a_k$, $b=b_k$ and the following associations:\\
  $\bullet$ for $\mu\in(0,\half)$, $\nu\in(0,\third)$:
    \[
    \theta_a = \bigfrac{2\xi L_1}{\tau\vartheta_L^2}(1-2\mu),
    \ms
    \theta_b = \bigfrac{2\xi L_2}{3\tau\vartheta_Q^3(1-3\nu)},
    \ms
    \theta =  \bigfrac{4\xi (f(x_0)-\flow)}{\tau};
    \]
  $\bullet$ for $\mu=\half$, $\nu\in(0,\third)$:
    \[
    \theta_a = \bigfrac{2\xi L_1}{\tau\vartheta_L^2},
    \ms
    \theta_b = \bigfrac{2\xi L_2}{3\vartheta_Q^3\tau(1-3\nu)},
    \ms
    \theta =  \bigfrac{4\xi (f(x_0)-\flow)}{\tau};
    \]
  $\bullet$ for $\mu\in(\half,1)$, $\nu\in(0,\third)$:
    \[
    \theta_a = 1,
    \ms
    \theta_b = \bigfrac{2\xi L_2}{3\tau\vartheta_Q^3(1-3\nu)},
    \ms
    \theta = \bigfrac{4\xi (f(x_0)-\flow)}{\tau}
    + \bigfrac{2\xi L_1}{\tau\vartheta_L^2}\cdot\frac{\varsigma^{1-2\mu}}{2\mu-1};
    \]
  $\bullet$ for $\mu\in(0,\half)$, $\nu=\third$:
    \[
    \theta_a = \bigfrac{2\xi L_1}{\tau\vartheta_L^2(1-2\mu)},
    \ms
    \theta_b = \bigfrac{2\xi L_2}{3\tau\vartheta_Q^3},
    \ms
    \theta =  \bigfrac{4\xi (f(x_0)-\flow)}{\tau};
    \]
  $\bullet$ for $\mu=\half$, $\nu=\third$:
    \[
    \theta_a = \bigfrac{2\xi L_1}{\tau\vartheta_L^2},
    \ms
    \theta_b = \bigfrac{2\xi L_2}{3\tau\vartheta_Q^3},
    \ms
    \theta =  \bigfrac{4\xi (f(x_0)-\flow)}{\tau};
    \]
  $\bullet$ for $\mu\in(\half,1)$, $\nu=\third$:
    \[
    \theta_a = \bigfrac{2\xi L_1}{\tau\vartheta_L^2},
    \ms
    \theta_b = \bigfrac{2\xi L_2}{3\vartheta_Q^3\tau};
    \ms
    \theta =  \bigfrac{4\xi (f(x_0)-\flow)}{\tau}
       +\bigfrac{2\xi L_1}{\tau\vartheta_L^2}\cdot\frac{\varsigma^{1-2\mu}}{2\mu-1};
    \]
  $\bullet$ for $\mu\in(0,\half)$, $\nu\in (\third,1)$:
    \[
    \theta_a = \bigfrac{2\xi L_1}{\tau\vartheta_L^2(1-2\mu)},
    \ms
    \theta_b = 1,
    \ms
    \theta =  \bigfrac{4\xi (f(x_0)-\flow)}{\tau}
    + \bigfrac{2\xi L_2}{3\tau\vartheta_Q^3}\cdot\frac{\varsigma^{1-3\nu}}{3\nu-1};
    \]
  $\bullet$ for $\mu=\half$, $\nu\in (\third, 1)$:
    \[
    \theta_a = \bigfrac{2\xi L_1}{\tau\vartheta_L^2},
    \ms
    \theta_b = 1,
    \ms
    \theta =  \bigfrac{4\xi (f(x_0)-\flow)}{\tau}
    + \bigfrac{2\xi L_2}{3\tau\vartheta_Q^3}\cdot\frac{\varsigma^{1-3\nu}}{3\nu-1};
    \]
  $\bullet$ for $\mu\in(\half,1)$, $\nu\in (\third,1)$:
    \[
    \theta_a = 1,
    \ms
    \theta_b = 1,
    \ms
    \theta =  \bigfrac{4\xi (f(x_0)-\flow)}{\tau}
       + \bigfrac{2\xi L_1}{\tau\vartheta_L^2}\cdot\frac{\varsigma^{1-2\mu}}{2\mu-1}
       + \bigfrac{2\xi L_2}{3\tau\vartheta_Q^3}\cdot\frac{\varsigma^{1-3\nu}}{3\nu-1}.
       \]
   As a consequence of applying Lemma~\ref{techsum}, we obtain that there
   exists positive constants\footnote{We choose them to be at least one, in
     order to cover the cases where $a_k\leq 1$ or $b_k\leq 1$ mentioned
     before \req{arglog2} and \req{arglog3}.}
   $\kap{1rst}\geq 1$ and $\kap{2nd}\geq 1 $ only depending on
   problem-related quantities and on $\varsigma$, $\xi$, $\mu$ and $\nu$ such
   that, for all $k\geq 0$,
   \beqn{partial-avgs}
   a_k \leq \kap{1rst}
   \tim{ and }
   b_k \leq \kap{2nd}.
   \eeqn
   We also have, from the mechanism of Step~3 of the algorithm (see
   \req{LQ-cond}) and \req{akbk-def}, that
   \[
   \sum_{j=0}^k \|g_j\|^2
   = \sum_{\mystack{j=0}{j\in\calK^L}}^k\|g_j\|^2 + \sum_{\mystack{j=0}{j\in\calK^Q}}^k\|g_j\|^2
   \leq
   \sum_{\mystack{j=0}{j\in\calK^L}}^k\|g_j\|^2+\sum_{\mystack{j=0}{j\in\calK^Q}}^k\hatphi_j^3
   \leq \half(a_k+b_k)
   \leq \half(\kap{1rst} + \kap{2nd})
   \]
   and
   \[
   \sum_{j=0}^k \hatphi_j^3
   = \sum_{\mystack{j=0}{j\in\calK^L}}^k\hatphi_j^3 + \sum_{\mystack{j=0}{j\in\calK^Q}}^k\hatphi_j^3
   \leq \sum_{\mystack{j=0}{j\in\calK^L}}^k\|g_j\|^2+\sum_{\mystack{j=0}{j\in\calK^Q}}^k\hatphi_j^3
   \leq \half(a_k+b_k)
   \leq \half(\kap{1rst} + \kap{2nd}).
   \]
   These two inequalities in turn imply that, for all $k\geq 0$,
   \[
   (k+1) \average_{j\in\iiz{k}} \|g_j\|^2 \leq  \half(\kap{1rst} + \kap{2nd})
   \tim{ and }
   (k+1) \average_{j\in\iiz{k}} \hatphi_j^3 \leq \half(\kap{1rst} + \kap{2nd}),
   \]
   and the desired results follow with $\kap{ASTR2} = \half(\kap{1rst} + \kap{2nd})$.
} %epr

Comments:
\begin{enumerate}
\item Note that $\hatphi_k < \phi_k$ only when $\phi_k > \xi$. Thus, if
  $\phi_k$ is bounded\footnote{Which is the case if $\|g_k\|\leq \kappa_g$ (as we
  will require in Section~\ref{ds-s}) since then
  $
  \phi_k \leq \|g_k\| + \half \|H_k\|
  \leq \kappa_g + \half L_1.$}, one can choose $\xi$
  large enough to ensure that $\phi_k=\hatphi_k$ for all $k$, and therefore
  that $\min_{j\iiz{k}}\phi_j \leq \kap{ASTR2}/(k+1)^\third$. In practice,
  $\xi$ can be used to tune the algorithm's sensitivity to second-order
  information.
\item If the $k$-orders of convergence specified by \req{min-orders} are translated
  in $\epsilon$-orders, that is numbers of iterations/evaluations to achieve
  $\|g(x_k\|\leq \epsilon_1$ and $\phi_k = \hatphi_k \leq \epsilon_2$, where
  $\epsilon_1$ and $\epsilon_2$ are precribed accuracies, we verify that at
  most $\calO(\epsilon_1^{-2})$ of them are needed to achieve the first of these
  conditions, while at most $\calO(\epsilon_2^{-3})$ are needed to achieve the
  second. As a consequence, at most
  $\calO(\max[\epsilon_1^{-2},\epsilon_2^{-3}])$ iterations/evaluations are needed to
  satisfy both conditions.  These orders are identical to the sharp bounds 
  known for the familiar trust-region methods (see
  \cite{GratSartToin08,CartGoulToin12d} or \cite[Theorems~2.3.7 and
    3.2.6]{CartGoulToin22}\footnote{This second of these theorems quotes an
    $\calO(\max[\epsilon_1^{-2}\epsilon_2^{-1},\epsilon_2^{-3}])$ order bound
    known for standard trust-region methods using first and second
    derivatives.}), or, for second-order optimality\footnote{Adaptive
    Regularization algorithms are faster for finding first-order points, as
    they find such points in $\calO(\epsilon_1^{-3/2})$ evaluations of the
    objective function and its gradient \cite{NestPoly06}, \cite[Theorem~3.3.9]{CartGoulToin22}.},
  for the Adaptive Regularization method (see
  \cite{NestPoly06}, \cite[Theorem~3.3.2]{CartGoulToin22}). This is quite remarkable because function
  values are essential in these two latter classes of algorithms to enforce
  descent, itself a crucial ingredient of existing convergence proofs.
\item While \req{min-orders} is adequate to allow a meaningful comparison of
  the global convergence rates with standard algorithms, as we just
  discussed, we note that \req{avg-orders} is stronger, because
  the average is of course a majorant of the minimum. One is then led to
  the question of whether such bounds in average can be proved for trust-region
  or adaptive regularization methods. As long as they haven't, \emph{the result
  presented here for second-order optimality can be viewed as one of the
  strongest available across all known methods using first and second derivatives.}
\item The expression of the constants is very intricate.  However it is
  remarkable that they do not explicitly depend on the problem dimension.
  However, and although a good sign, this does
  not tell the whole story and caution remains advisable, because the Lipschitz
  constants $L_1$ and $L_2$ may themselves hide this (potentially severe)
  dependence.
\item It is also remarkable that the bounds \req{avg-orders} and
  \req{min-orders} specify the same order of global convergence irrespective
  of the values of $\mu$ and $\nu$ in $(0,1)$, although these values do affect
  the constants involved.
\item The condition \req{LQ-cond} determining the choice of a linear (in
  $\calK^L$) or quadratic (in $\calK^Q$) step is only used at the very end of
  the theorem's proof, after \req{partial-avgs} has already been obtained.
  This means that other choice mechanisms are possible without affecting this
  last conclusion, which is enough to derive bounds on $\|g_j\|^2$ and
  $\hatphi_j^3$ averaged on iterations in $\calK^L$ and $\calK^Q$, respectively (rather than
  on all iterations).
\end{enumerate}

We now show that the bound \req{min-orders} is essentially sharp (in the sense
of \cite{CartGoulToin18a}, meaning that a lower bound on evaluation complexity exists which
is arbitrarily close to its upper bound) by following ideas of
\cite[Theorem~2.2.3]{CartGoulToin22} in an argument parallel to
that used in \cite{GratJeraToin22b} for the first-order bound.

\lthm{sharp2}{The bound \req{min-orders} is essentially sharp in that, for each
$\mu,\nu\in(0,1)$, $\vartheta_L=\vartheta_Q=1$ and each
  $\varepsilon\in(0,\sfrac{2}{3})$, there exists a univariate function
  $f_{\mu,\nu,\varepsilon}$ satisfying AS.1--AS.3 such that, when applied to
  minimize $f_{\mu,\nu,\varepsilon}$ from the origin,  the \al{ASTR2}
  algorithm  with new{\req{w-adag-L}-\req{w-adag-Q}} produces second-order
  optimality measures given by 
\beqn{is-sharp2}
\phi_k = \hatphi_k = \!\!\min_{j\in\iiz{k}}\hatphi_j = \frac{1}{(k+1)^{\third+\varepsilon}}.
\eeqn
}

\proof{We start by constructing $\{x_k\}$ for which
$f_{\mu,\nu,\varepsilon}(x_k) = f_k$, $\nabla_x^1f_{\mu,\nu,\varepsilon}(x_k) = g_k$ and
$\nabla_x^2 f_{\mu,\nu,\varepsilon}(x_k) = H_k$ for associated sequences of function, gradient and
Hessian values $\{f_k\}$, $\{g_k\}$ and $\{H_k\}$, and then apply Hermite
interpolation to exhibit the function $f_{\mu,\nu,\varepsilon}$ itself.  We
select an arbitrary $\varsigma > 0$ and define, for
$k\geq 0$,
\beqn{gkHk-def}
g_k \eqdef 0,
\tim{ and }
H_k = -\frac{2}{(k+1)^{\third+\varepsilon}},
\eeqn
from which we deduce, using \req{phi-def}, that, for $k> 0$,
\[
\phi_k=\hatphi_k= \frac{1}{(k+1)^{\third+\varepsilon}}.
\]
Since $\phi_k^3 > 0 = \|g_k\|^2$, we set
\beqn{sk-def}
s_k = s_k^Q \eqdef \frac{1}{(k+1)^{\third+\varepsilon}[\varsigma+\sum_{j=0}^k \hatphi_j^3]^\nu}, %\ms(k>0)
\eeqn
which is the exact minimizer of the quadratic model within the trust region,
yielding that, for $k\geq 0$,
\beqn{gksk}
\Delta q_k \eqdef
\left| g_ks_k + \half H_ks_k^2 \right|
= \frac{1}{(k+1)^{1+3\varepsilon}\big[\varsigma+\sum_{j=0}^k \hatphi_j^3)\big]^{2\nu}}
\leq \frac{1}{(k+1)^{1+3\varepsilon}},
\eeqn
where we used the fact that $\varsigma+\sum_{j=0}^k \hatphi_j^3>
\varsigma+\hatphi_0>1$ to deduce the last inequality.
We then define, for all $k\geq0$,
\beqn{xk-def}
x_0 = 0,
\ms
x_{k+1} = x_k+s_k \ms(k\geq0)
\eeqn
and
\beqn{fk-def}
f_0 = \zeta(1+3\varepsilon),
\ms
f_{k+1} = f_k - \Delta q_k \ms (k\geq0),
\eeqn
where $\zeta(\cdot)$ is the Riemann zeta function.
Observe that the sequence $\{f_k\}$ is decreasing and that, for all $k\geq 0$,
\beqn{fk-sum}
f_{k+1}
= f_0 - \bigsum_{k=0}^k\Delta q_k
\geq f_0 - \bigsum_{k=0}^k\frac{1}{(k+1)^{1+3\varepsilon}}
\geq f_0 - \zeta(1+3\varepsilon),
\eeqn
where we used \req{fk-def} and \req{gksk}. Hence
\req{fk-def} implies that
\beqn{fk-bound}
f_k \in [0, f_0] \tim{for all} k\geq 0.
\eeqn
Also note that, using \req{fk-def},
\beqn{dfok}
|f_{k+1} - f_k + \Delta q_k| = 0,
\eeqn
while, using \req{gkHk-def},
\beqn{dgok}
|g_{k+1}-g_k| = 0
\ms (k \geq 0).
\eeqn
Moreover, using the fact that $1/x^{\third+\nu}$ is a convex function of $x$ over
$[1,+\infty)$, and that from \req{sk-def}
$s_k \geq \frac{1}{(k+1)^{\third+\nu} \left(\varsigma + k+1\right)^\nu}$,  we
derive that, for $k\geq0$,
\begin{align*}
|H_{k+1} - H_k| &= 2\left| \frac{1}{(k+2)^{\third+\nu}} - \frac{1}{(k+1)^{\third+\nu}} \right| \\
&\leq 2\left(\frac{1}{3}+\nu\right) \frac{1}{(k+1)^{\sfrac{4}{3} + \nu}} \\
&\leq \frac{8}{3} \, \frac{(\varsigma + k+1 )^\nu}{ (k+1)  (k+1)^{\third + \nu} (\varsigma + k+1)^\nu} \\
&\leq \frac{8}{3} \, \frac{(\varsigma + k+1 )^\nu}{k+1 } s_k \\
&\leq \frac{8}{3} \, \left(\varsigma + 2 \right)^\nu s_k .
\end{align*}
These last bounds with \req{fk-bound}, \req{dfok} and \req{dgok} allow us to use
standard Hermite interpolation on the data given by $\{f_k\}$, $\{g_k\}$ and
$\{H_k\}$: see, for instance, Theorem~A.9.1 in \cite{CartGoulToin22} with $p=2$ and
\[
\kappa_f = \max\left[\frac{8}{3}(\varsigma+2)^\nu, f_0,2\right]
\]
(the second term in the max bounding $|f_k|$ because of
\req{fk-bound} and the third bounding $|g_k|$ and $|H_k|$ because of \req{gkHk-def}).
We then deduce that there exists a twice continuously differentiable function $f_{\mu,\nu,\varepsilon}$
from $\Re$ to $\Re$ with Lipschitz continuous gradient and Hessian (i.e. satisfying
AS.1 and AS.2) such that, for $k\geq 0$,
\[
f_{\mu,\nu,\varepsilon}(x_k) = f_k,\ms \nabla_x^1f_{\mu,\nu,\varepsilon}(x_k) = g_k
\tim{ and } \nabla_x^2f_{\mu,\nu,\varepsilon}(x_k) = H_k.
\]
Moreover, the range of $f_{\mu,\nu,\varepsilon}$ is constant independent
of $\varepsilon$, hence guaranteeing AS.3.
The definitions \req{gkHk-def}, \req{sk-def}, \req{xk-def} and \req{fk-def}
imply that the sequences $\{x_k\}$,  $\{f_k\}$,  $\{g_k\}$ and $\{H_k\}$ can be seen as generated by the
\al{ASTR2} algorithm applied to $f_{\mu,\nu,\varepsilon}$, starting from $x_0=0$.
} %epr

\begin{figure}[htb] % produced by slowexample.m
\centerline{
\includegraphics[height=5cm,width=5.2cm]{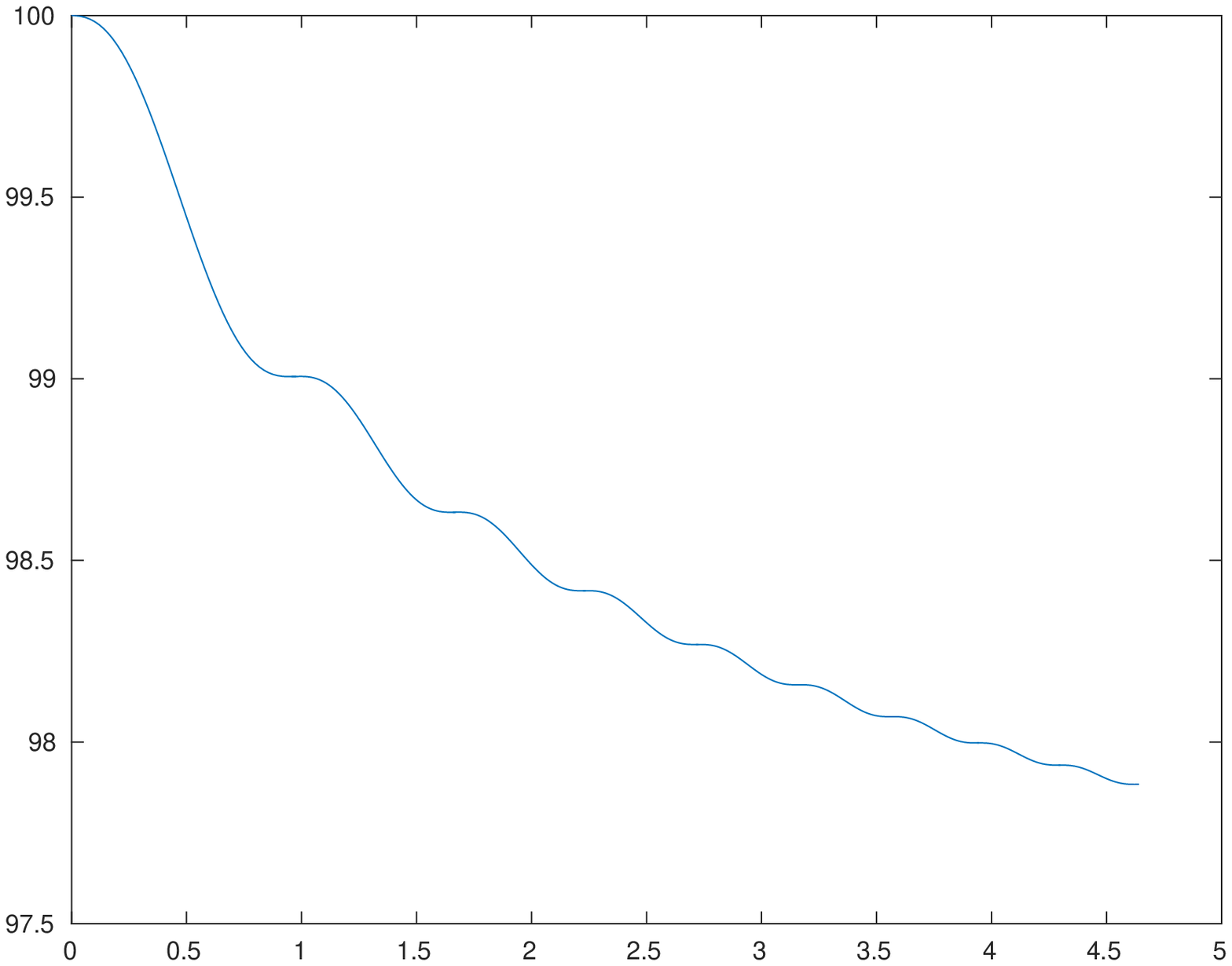}
\includegraphics[height=5cm,width=5.2cm]{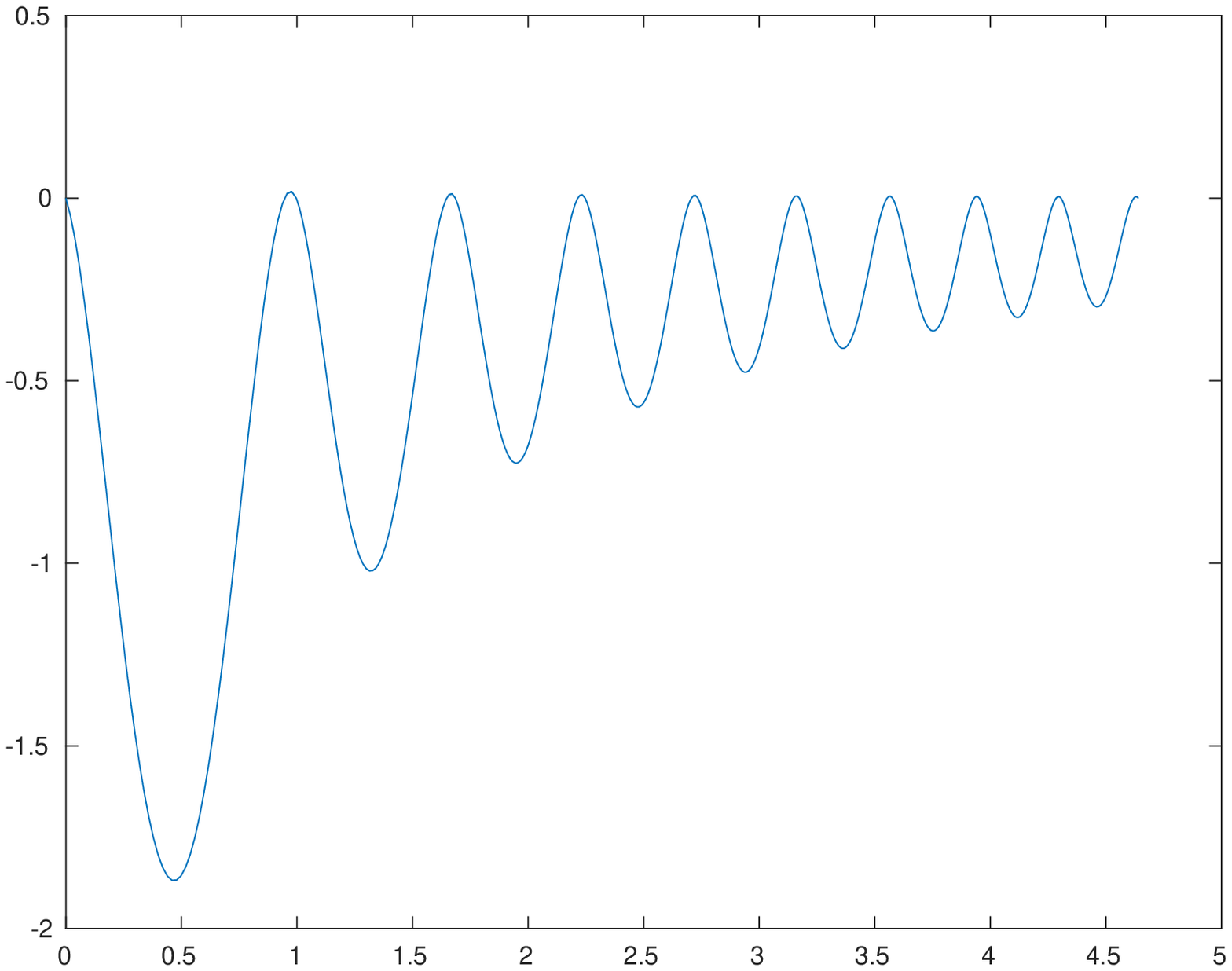}
\includegraphics[height=5cm,width=5.2cm]{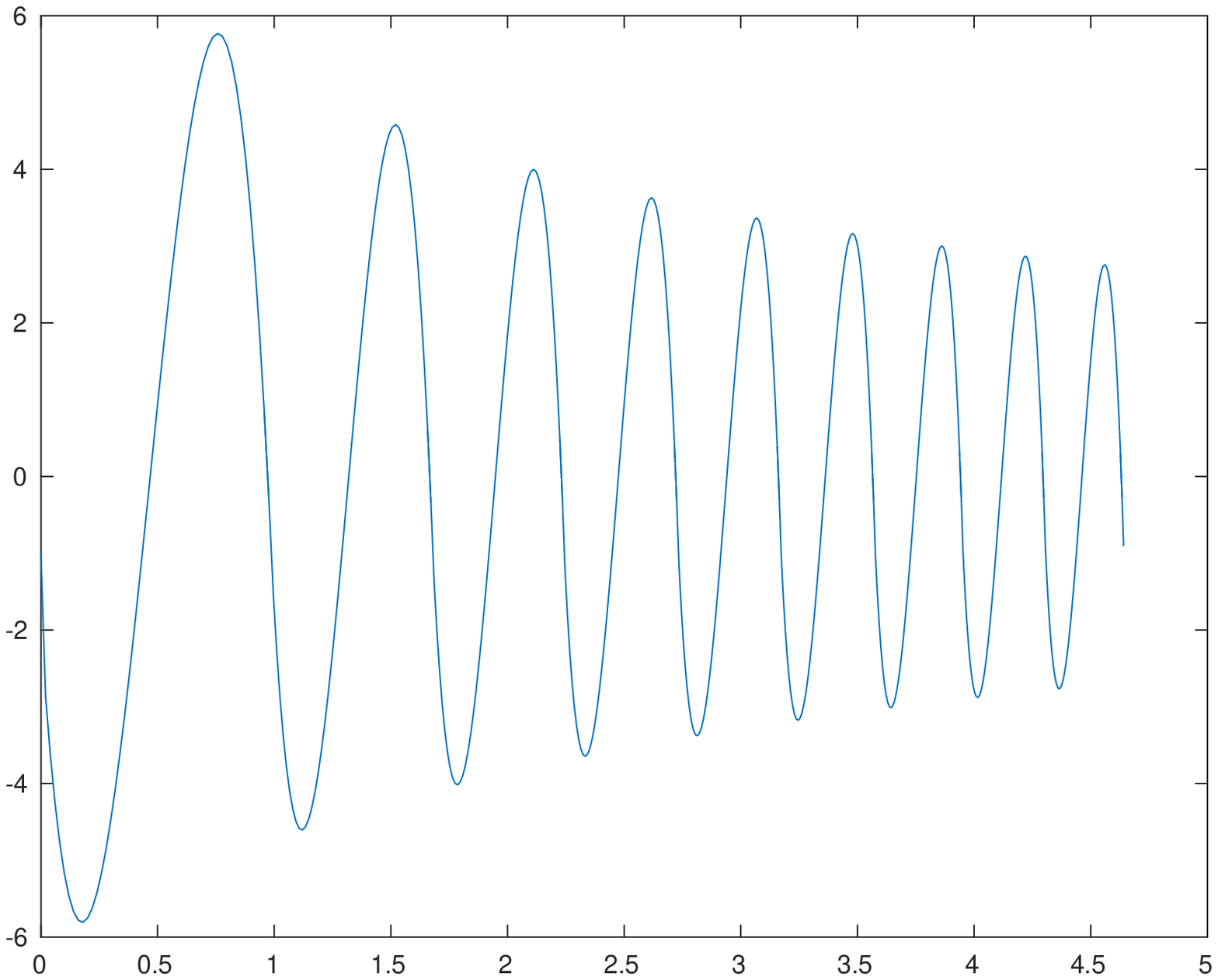}
}
\caption{\label{figure:slowex}
  The function $f_{\mu,\nu,\varepsilon}(x)$ (left), its gradient $\nabla_x^1f_{\mu,\nu,\varepsilon}(x)$
  (middle) and its Hessian $\nabla_x^2f_{\mu,\nu,\varepsilon}(x)$ (right) plotted as a
  function of $x$, for the first 10 iterations of
  the \al{ASTR2} algorithm with \req{w-adag-L}-\req{w-adag-Q} ($\mu = \sfrac{1}{2}$, $\nu=\third$,
  $\varepsilon =\varsigma = \sfrac{1}{100}$, $\vartheta_L=\vartheta_Q=1$)}   
\end{figure}

\noindent
Figure~\ref{figure:slowex} shows the behaviour of $f_{\mu,\nu,\varepsilon}(x)$ for $\mu =
\half$, $\nu = \third$, $\vartheta_L=\vartheta_Q=1$ and $\varepsilon = \varsigma=\sfrac{1}{100}$, its gradient
and Hessian, as resulting from the first 10 iterations of the \al{ASTR2}
algorithm with \req{w-adag-L}-\req{w-adag-Q}. (We have chosen to shift $f_0$ to 100 in order to
avoid large numbers on the vertical axis of the left panel.) Due to the slow
convergence of the series $\sum_j 1/j^{\frac{1}{1+3/100}}$, illustrating the
boundeness of $f_0-f_{k+1}$ would require many more iterations. One also notes
that the gradient is not monotonically increasing, which implies that
$f_{\mu,\nu,\varepsilon}(x)$ is nonconvex, as can be verified in the left
panel. Note that the unidimensional nature of the example
is not restrictive, since it is always possible to make the value of its
objective function and gradient independent of all dimensions but one.
Also note that, as was the case in \cite{GratJeraToin22b}, the argument of
Theorem~\ref{sharp2} fails for $\varepsilon=0$ since then the sums in \req{fk-sum}
diverge when $k$ tends to infinity.

Note that, because
\begin{align*}
\sum_{j=0}^k \frac{1}{(j+1)^{\third+\varepsilon}}
&\geq \int_0^k \frac{dj}{(j+2)^{\third+\varepsilon}}\\
&=\frac{3}{2+3\varepsilon}\left[\frac{k+2}{(k+2)^{\third+\varepsilon}}-2\right]\\
&\geq \frac{3}{2(2+3\varepsilon)}\left[\frac{k+1}{(k+1)^{\third+\varepsilon}}-2\right],
\end{align*}
one deduces that
\[
\average_{j\in\iiz{k}}\hatphi_j \geq \frac{3}{2(2+3\varepsilon)}\left[\frac{1}{(k+1)^{\third+\varepsilon}}-\frac{2}{k+1}\right],
\]
which, when compared to \req{min-orders}, reflects the (slight) difference in
strength between \req{avg-orders} and \req{min-orders}.

\numsection{A ``divergent stepsize'' \tal{ASTR2} subclass}\label{ds-s}

A ``divergent stepsize'' first-order method was analyzed in \cite{GratJeraToin22b},
motivated by its good practical behaviour in the stochastic context
\cite{GratJeraToin22a}. For coherence, we now present and analyze a similar
variant, this time for second-order optimality. This requires the following
additional assumption.
\begin{description}
\item[AS.4:] there exists a constant $\kappa_g >0$ such that, for all $x$,
   $\|g(x)\|_\infty\le \kappa_g$.
\end{description}

\lthm{theorem:div-series}{ Suppose that AS.1--AS.3 and AS.4 hold and that the
\al{ASTR2} algorithm is applied to problem \req{problem}, where,
the scaling factors $w_{i,k}$ are chosen such that, for some power
parameters $0< \nu_1 \leq \mu_1 < 1$ and $0<\nu_2\leq \mu_ 2<\half$, some constants
$\varsigma\in (0,1]$ and
$\kappa_w\geq \max[1,\varsigma]$,  all $i\in\ii{n}$ and all $k\geq0$,
\beqn{wk-divs}
0 < \varsigma\, (k+1)^{\nu_1} \le w^L_{k} \le \kappa_w \,(k+1)^{\mu_1}
\tim{ and }
0 < \varsigma\, (k+1)^{\nu_2} \le w^Q_{k} \le \kappa_w \,(k+1)^{\mu_2}.
\eeqn
Let $\psi_k\eqdef \min[1,\max[\|g_k\|^2,\phi_k^3]]$.
Then, for any $\theta \in (0,\quarter \tau)$ and $k >j_\theta$,
\beqn{psi-bound-divs}
\min_{j\in\iibe{j_\theta}{k}}\psi_k
\le \kappa_\diamond(\theta) \bigfrac{(k+1)^{\max[\mu_1,2\mu_2]}}{k-j_\theta}
\le \bigfrac{\kappa_\diamond(\theta)(j_\theta+1)}{(k+1)^{1-\max[\mu_1,2\mu_2]}},
\eeqn
where
\beqn{jtheta-divs}
j_\theta \eqdef
\max\left[
    \left(\frac{L_1}{2\varsigma(1-\theta)}\right)^{\sfrac{1}{\nu_1}},
    \left(\frac{2(1+L_1)}{\varsigma}\right)^{\sfrac{1}{\nu_2}},
    \left(\frac{L_2}{3\varsigma(\quarter\tau-\theta)}\right)^{\sfrac{1}{\nu_2}},
    \left(\frac{L_2\xi}{3\varsigma^2(\quarter\tau-\theta)}\right)^{\sfrac{1}{2\nu_2}}
    \right]
\eeqn
and
\[
\kappa_\diamond(\theta)
\eqdef \left\{ \frac{\kappa_w^2}{\theta}\left( f(x_0)-\flow +
      (j_\theta+1)\max\left[\frac{L_1\kappa_g^2}{2\varsigma^2},\frac{L_2\xi^3}{3\varsigma^3}\right] \right)\right\}^\third.
\]
}

\proof{Consider an arbitrary $\theta \in (0,\quarter \tau)$ and note that
  AS.4, \req{wk-divs} and the definition of $\hatphi_k$ imply that
  \beqn{w-bounds}
  w^L_k \in [\varsigma^\mu, \kappa_w(k+1)^{\mu_1}]
  \tim{ and }
  w^Q_k \in [\varsigma^\nu, \kappa_w(k+1)^{\mu_2} ].
  \eeqn
  If we define $j_\theta$ by \req{jtheta-divs}, we immediately obtain from AS.4
  and Lemma~\ref{lemma:fdecr} (where we neglect the first term in the right-hand
  sides of \req{ffdecr0} and \req{ffdecr1}) that
  \beqn{fover}
  f(x_{j_\theta+1}) \leq f(x_0) + (j_\theta+1)\kap{over}
  \tim{ where }
  \kap{over} =
  \max\left[\frac{L_1\kappa_g^2}{2\varsigma^2},\frac{L_2\xi^3}{3\varsigma^3}\right].
  \eeqn
  If we choose $j > j_\theta$, one then verifies that the definition of
  $j_\theta$ in \req{jtheta-divs}, the bounds \req{ffdecr0} and \req{ffdecr1}
  and the definition \req{wk-divs} together ensure that 
  \[
  f(x_{j+1}) -f(x_j)
  \leq \left\{ \begin{array}{ll}
  -\theta\bigfrac{\|g_k\|^2}{w^L_k}    & \tim{if } j\in\calK^L,\\*[3ex]
  -\theta\bigfrac{\phi_k^3}{(w^Q_k)^2} & \tim{if } j\in\calK^Q.
  \end{array}\right.
  \]
  Using now the mechanism of Step~3, the
  definition of $\psi_k$, \req{wk-divs} and the inequality
  $\kappa_w\geq 1$, we obtain that, for $j > j_\theta$
  \beqn{fdecr-large}
  f(x_j)-f(x_{j+1})
  \geq \theta\psi_j\min\left[\frac{1}{w^L_k},\frac{1}{(w^Q_k)^2}\right]
  \geq \frac{\theta\psi_j}{\kappa_w^2(j+1)^{\max[\mu_1,2\mu_2]}}.
  \eeqn
  As a consequence, we obtain from \req{fover} and the summation of \req{fdecr-large}
  for $j\in \iibe{j_\theta+1}{k}$ that, for $k > j_\theta$,
  \[
  f(x_0)-f(x_{j+1})
  \geq -(j_\theta+1)\kap{over}
       +\sum_{j=j_\theta+1}^k\frac{\theta\psi_j}{\kappa_w^2(j+1)^{\max[\mu_1,2\mu_2]}}.
  \]
  We therefore deduce, using AS.3, that
  \[
  (k-j_\theta)\min_{j_{\theta,\max}+1,\ldots,k}\psi_j
  \leq \sum_{j=j_\theta+1}^k\psi_j
  \leq \frac{\kappa_w^2(k+1)^{\max[\mu_1,2\mu_2]}}{\theta}\Big[ f(x_0)-\flow +
      (j_\theta+1)\kap{over} \Big],
  \]
  and \req{psi-bound-divs} follows.
} % epr

\noindent
This theorem gives a bound on the rate at which the combined optimality
measure $\psi_k$ tends to zero, and this bound is slightly worse than but
close to what we obtained in the previous section whenever
$\max[\mu_1,2\mu_2]$ approaches zero.

Using the methodology of Theorem~\ref{sharp2}, we now show that the bound
\req{psi-bound-divs} is also essentially sharp. 

\lthm{sharp2-divs}{The bound \req{psi-bound-divs} is essentially sharp in that, for each
  $\mu = (\mu_1,\mu_2)$, each $\nu = (\nu_1,\nu_2)$ with
  $0<\nu_1\leq\mu_1<1$ and $0<\nu_2\leq\mu_2<\half$
  and each $\varepsilon \in (0,1-\third(1-2\mu_2))$, there exists a univariate function $h_{\mu,\nu,\varepsilon}$ 
satisfying AS.1--AS.4 such that, when applied to minimize $h_{\mu,\nu,\varepsilon}$ from the
origin, the \al{ASTR2} algorithm  with \req{wk-divs} produces 
second-order optimality measures given by
\beqn{is-sharp2-divs}
\phi_k = \hatphi_k = \psi_k = \min_{j\in\iiz{k}}\psi_j = \frac{1}{(k+1)^{\third(1-2\mu_2)+\varepsilon}}.
\eeqn
}

\proof{As above, we start by defining, for $k\geq 0$,
$\gamma = \third(1-2\mu_2)+\varepsilon$, $w_k= \kappa_w(k+1)^{\mu_2}$,
and, for $k\geq 0$,
\beqn{gkHk-def-divs}
g_k \eqdef 0,
\tim{ and }
H_k = -\frac{2}{(k+1)^\gamma},
\eeqn
which then implies, using \req{phi-def} that, for $k> 0$,
\beqn{hpk}
\phi_k=\hatphi_k= \frac{1}{(k+1)^\gamma}.
\eeqn
Given these definitions and because $\hatphi^3_k>0 = \|g_k\|^2$, we set
\beqn{sk-def-divs}
s_k = s_k^Q
\eqdef \frac{1}{(k+1)^\gamma[\kappa_w(k+1)^{\mu_2}]}
= \frac{1}{\kappa_w(k+1)^{\gamma+\mu_2}},
\eeqn
yielding that, for $k>0$,
\beqn{gksk-divs}
\Delta q_0 \eqdef \frac{1}{(\varsigma+1)^{2\nu}}
\tim{ and }
\Delta q_k \eqdef
\left| g_ks_k + \half H_ks_k^2 \right|
= \frac{1}{\kappa_w^2(k+1)^{3\gamma+2\mu_2}}
\leq \frac{1}{(k+1)^{3\gamma+2\mu_2}},
\eeqn
where we used the fact that $\kappa_w \geq 1$ to deduce the last inequality.
We then define, for all $k\geq0$,
\beqn{xk-def-divs}
x_0 = 0,
\ms
x_{k+1} = x_k+s_k \ms(k>0)
\eeqn
and
\beqn{fk-def-divs}
h_0 = \zeta(3\gamma+2\mu_2)
\tim{ and }
h_{k+1} = h_k - \Delta q_k \ms (k \geq 0),
\eeqn
where $\zeta(\cdot)$ is the Riemann zeta function. Note that, since
$\gamma>1-2\mu_2$, the argument $3\gamma+2\mu_2$ of $\zeta$ is strictly larger
than one and $\zeta(3\gamma+2\mu_2)$ is finite.
Observe also that the sequence $\{h_k\}$ is decreasing and that, for all $k\geq 0$,
\beqn{fk-sum2}
h_{k+1}
= h_0 - \bigsum_{k=0}^k\Delta q_k
\geq h_0 - \bigsum_{k=0}^k\frac{1}{(k+1)^{3\gamma+2\mu_2}}
\geq h_0 - \zeta(3\gamma+2\mu_2),
\eeqn
where we used \req{fk-def} and \req{gksk}. Hence
\req{fk-def} implies that
\beqn{fk-bound-divs}
h_k \in [0, h_0] \tim{for all} k\geq 0.
\eeqn
Also note that, using \req{fk-def},
\beqn{dfok-divs}
|h_{k+1} - h_k + \Delta q_k| = 0,
\eeqn
while, using \req{gkHk-def},
\beqn{dgok-divs}
|g_{k+1}-g_k| = 0
\ms (k \geq 0).
\eeqn
Moreover, using the fact that $1/x^\gamma$ is a convex function of $x$ over
$[1,+\infty)$ and  \req{sk-def-divs},  we derive that, for $k\geq0$,
\[
|H_{k+1} - H_k| = 2\left| \frac{1}{(k+2)^\gamma} - \frac{1}{(k+1)^\gamma} \right| 
\leq \frac{2\gamma}{(k+1)^{1+\gamma}} 
\leq \frac{2\gamma\kappa_w(k+1)^{\mu_2}}{k+1} \,s_k
\leq 2 \gamma\kappa_w s_k.
\]
This bound with \req{fk-bound-divs}, \req{dfok-divs} and \req{dgok-divs} once
more allow us to use standard Hermite interpolation on the data given by $\{h_k\}$, $\{g_k\}$ and
$\{H_k\}$, as stated in \cite[Theorem~A.9.1]{CartGoulToin22} with $p=2$ and
\[
\kappa_f = \max\left[2\gamma\kappa_w, h_0,2\right]
\]
(the second term in the max bounds $|h_k|$ because of
\req{fk-bound-divs} and the third bounds both $|g_k|$ and $|H_k|$ because of \req{gkHk-def-divs}).
As a consequence, there exists a twice continuously differentiable function $h_{\mu,\nu,\varepsilon}$
from $\Re$ to $\Re$ with Lipschitz continuous gradient and Hessian (i.e. satisfying
AS.1 and AS.2) such that, for $k\geq 0$,
\[
h_{\mu,\nu,\varepsilon}(x_k) = h_k,\ms \nabla_x^1h_{\mu,\nu,\varepsilon}(x_k) = g_k
\tim{ and } \nabla_x^2h_{\mu,\nu,\varepsilon}(x_k) = H_k.
\]
Moreover, the ranges of $h_{\mu,\nu,\varepsilon}$ and its derivatives is constant independent
of $\gamma$, hence guaranteeing AS.3 and AS.4.
Thus \req{gkHk-def-divs}, \req{sk-def-divs}, \req{xk-def-divs} and \req{fk-def-divs}
imply that the sequences $\{x_k\}$,  $\{h_k\}$,  $\{g_k\}$ and $\{H_k\}$ can be seen as generated by the
\al{ASTR2} algorithm applied to $h_{\mu,\nu,\varepsilon}$, starting from $x_0=0$. The
first bound of \req{is-sharp2-divs} then results from \req{hpk} and the
definition of $\gamma$.
}%epr

\begin{figure}[htb] % produced by slowexample2.m
\centerline{
\includegraphics[height=5cm,width=5.2cm]{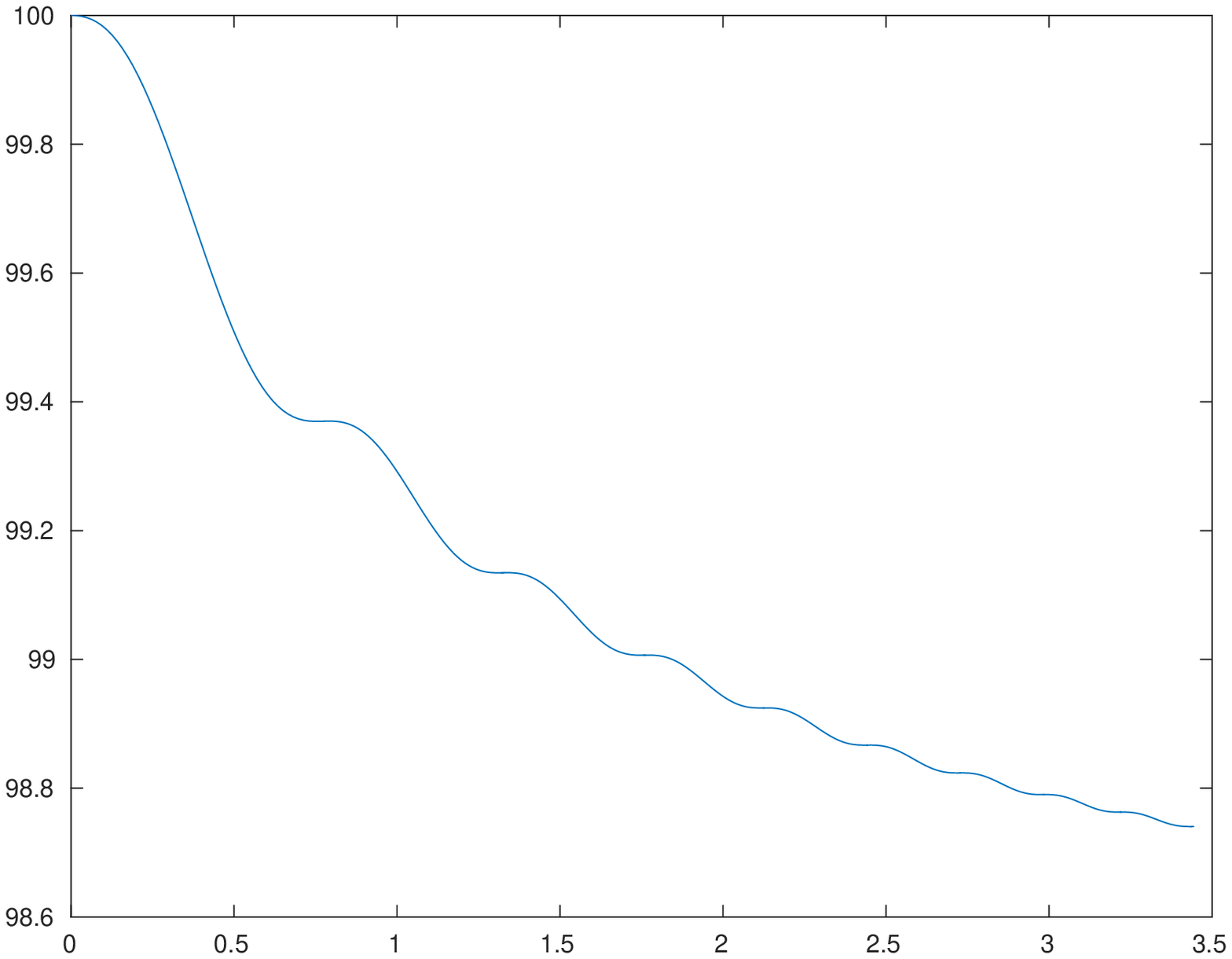}
\includegraphics[height=5cm,width=5.2cm]{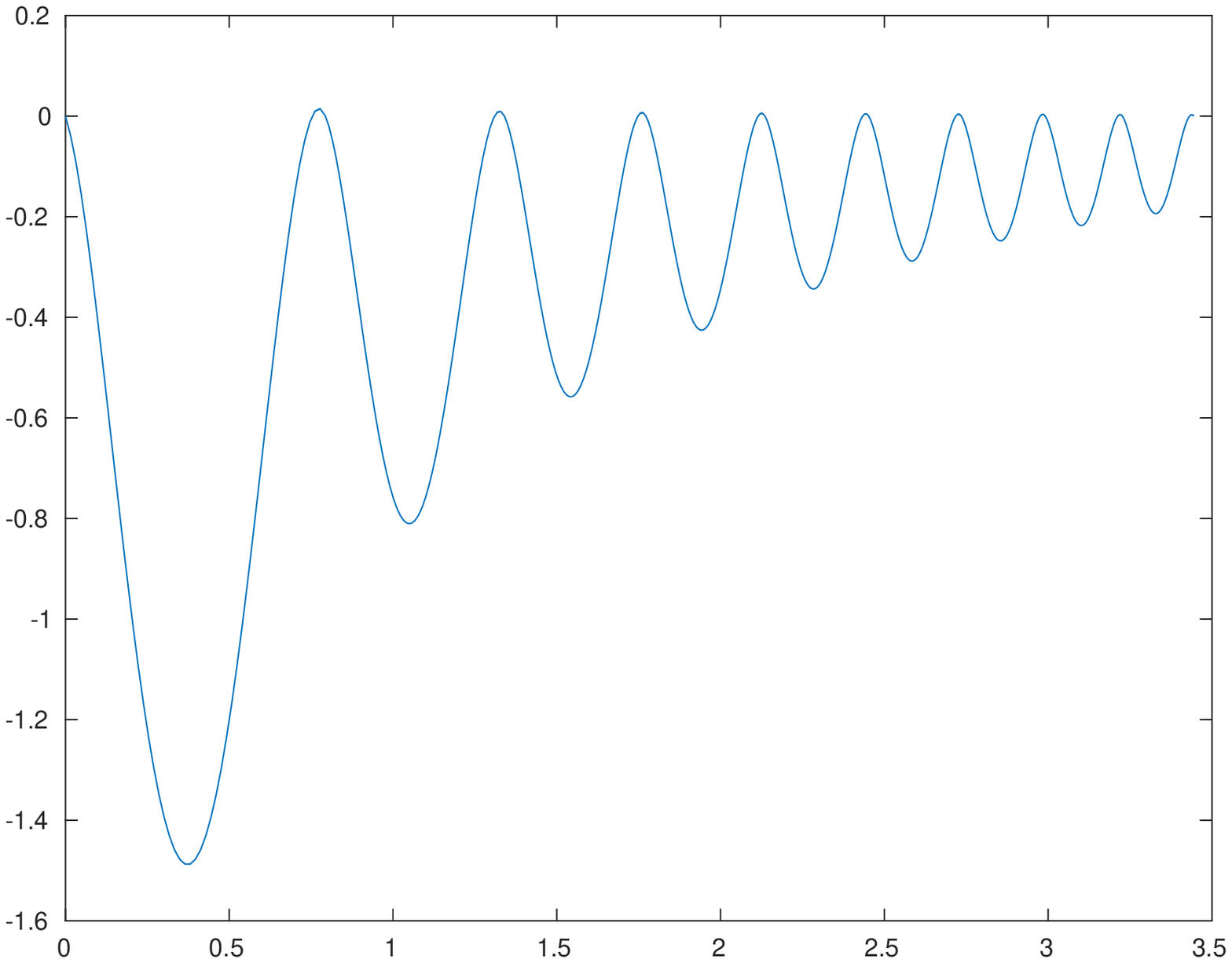}
\includegraphics[height=5cm,width=5.2cm]{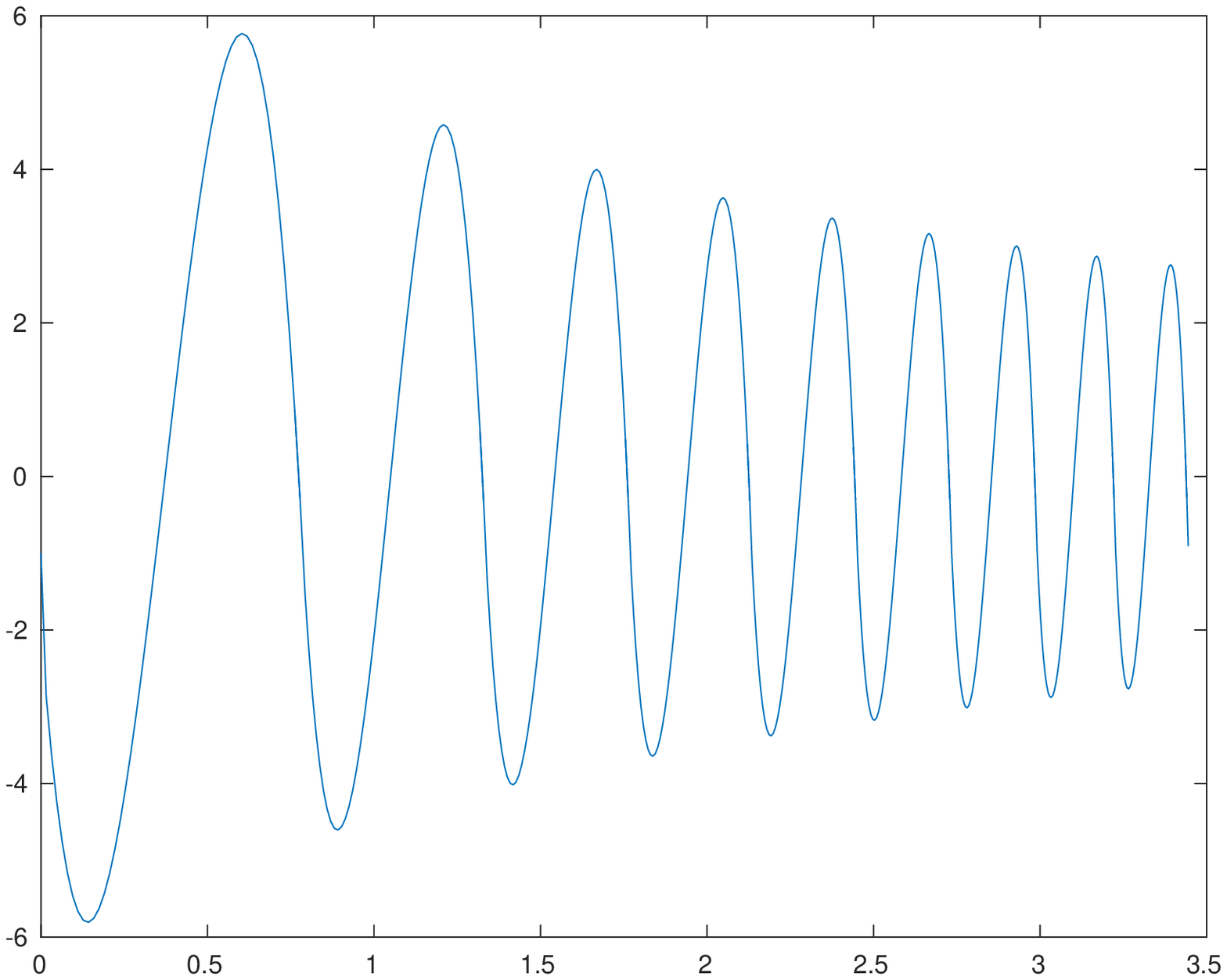}
}
\caption{\label{figure:slowex2}
  The function $h_{\mu,\nu,\varepsilon}(x)$ (left), its gradient $\nabla_x^1h_{\mu,\nu,\varepsilon}(x)$
  (middle) and its Hessian $\nabla_x^2h_{\mu,\nu,\varepsilon}(x)$ (right) plotted as a
  function of $x$, for the first 10 iterations of
  the \al{ASTR2} algorithm with \req{psi-bound-divs} ($\mu = \nu = (\sfrac{1}{2},\sfrac{1}{3})$)}
\end{figure}

\noindent
The behaviour of $h_{\mu,\nu,\varepsilon}$ is illustrated in
Figure~\ref{figure:slowex2}. It is qualitatively similar to that of
$f_{\mu,\nu,\varepsilon}$ shown in Figure~\ref{figure:slowex}, although the decrease
in objective-value is somewhat slower, as expected.
As in Section~\ref{adag-s}, note that the inequality
\[
\sum_{j=0}^k\frac{1}{(j+1)^\gamma}
\geq \int_0^k\frac{dj}{(j+2)^\gamma}
= \frac{1}{1-\gamma}\left[\frac{k+2}{(k+2)^\gamma}-1\right]
\geq \frac{1}{2(1-\gamma)}\left[\frac{k+1}{(k+1)^\gamma}-2\right]
\]
implies that
\[
\average_{j\in\iiz{k}}\psi_j \geq \frac{1}{2(1-\gamma)}\left[\frac{1}{(k+1)^{\third(1-2\mu_2)+\varepsilon}}-\frac{2}{k+1}\right],
\]
which has the same flavour as the second bound of \req{is-sharp2}.

\numsection{Second-order optimality in a subspace}\label{weak2nd-s}
  
While the \al{ASTR2} algorithms guarantee second-order optimality conditions,
they come at a computational price.  The key of this guarantee is of course
that significant negative curvature in any direction of $\Re^n$ must be
exploited, which requires evaluating the Hessian.  In addition, the optimality
measure $\phi_k$ and the step $s_k$ must also be computed.  However, these
computational costs may be judged excessive, so the question arises whether a potentially
cheaper algorithm is able to ensure a ``degraded'' or weaker form of
second-order optimality. Fortunately, the answer is positive: one can
guarantee second-order optimality in subspaces of $\Re^n$ at lower cost.

The first step is to assume that a subspace $\calS_k$ is of interest at
iteration $k$.  Then, instead of computing $\phi_k$ from \req{phik-2}, one can
choose to calculate
\[
\phi_k^{\calS_k}
= \max_{\mystack{\|d\|\le 1}{d\in\calS_k}} -\Big( g(x)^Td + \half d^TH(x)d \Big).
\]
Because the dimension of $\calS_k$ may be much smaller than $n$, the cost of
this computation may be significantly smaller than that of computing
$\phi_k$. The measure $\phi_k^{\calS_k}$ may for instance be obtained using a
Krylov-based method, as conjugate gradients \cite{HestStie52}, GLRT
\cite{GoulLuciRomaToin99} or variants thereof, where the minimum of the model
$T_{f,2}(x,d)$ within the trust region is derived iteratively in a sequence of
nested Krylov subspaces of increasing dimension, which tend to contain vector
along which curvature is extreme \cite[Chapter~9]{GoluvanL96}, thereby improving the quality of the
second-order guarantee compared to random subspaces. This process may then be
terminated before the subspaces fill $\Re^n$, should the calculation become too expensive or a
desired accuracy be reached. In addition, there is no need for $n_k$, the
dimension of the final Krylov space at iteration $k$ to be constant: it is
often kept very small when far from optimality. This technique has the added
benefit that the full Hessian is not evaluated, but only $n_k$
Hessian-times-vector products are needed, again significantly reducing the
computational burden. Calculating the step $s_k^Q$ for $k\in\calK^Q$ once
$\phi_k^{\calS_k}$ is known is also cheaper in a space of dimension $n_k$ much
less than $n$, especially since only a $\tau$-approximation is needed (see the
comments after the algorithm).

Importantly, the theory developped in the previous sections is
not affected by the transition from $\Re^n$ to $\calS_k$, except that now
the complexity bounds \req{avg-orders}-\req{min-orders} and
\req{psi-bound-divs} are no longer expressed using $\hatphi_k$ but now involve
$\hatphi_k^{\calS_k} = \min[ 1, \phi_k^{\calS_k}]$ instead. While clearly not as
powerful as the complete second-order guarantee in $\Re^n$, weaker
guarantees based on (Krylov) subspaces are often sufficient in practice and
make the \al{ASTR2} algorithm more affordable.
Note that, in the limit, one can even choose $\calS_k= \{0\}$ for all $k$, in
which case we can set $H_k=0$ for all $k$ and we do not obtain any second-order guarantee
(but the first-order complexity bounds remain valid, recovering results of \cite{GratJeraToin22b}).

\numsection{Conclusions}\label{concl-s}

We have introduced an OFFO algorithm whose global rate of convergence to
first-order minimizers is $\calO((k+1)^{-\half})$ while it
converges to second-order ones as $\calO((k+1)^{-\third})$. These
bounds are equivalent to the best known bounds for second-order optimality for
algorithms using objective-function evaluations, despite the latter exploiting
significantly more information. Thus we conclude that, from the point of
view of evaluation complexity at least, evaluating values of the objective
function is an unnecessary effort for efficiently finding second-order minimizers.  We have also discussed
another closely related algorithm, whose global rates of convergence can be
nearly as good. We have finally considered how weaker second-order guarantees
may be obtained at a much reduced computational cost.

We expect that extending our proposal to convexly constrained cases (for
instance to problems involving bounds on the variables) should be possible. As
in \cite[Chapter~12]{ConnGoulToin00}, the idea would be to restrict the model
minimization at each iteration to the intersection of the trust region with
the feasible domain, but this should of course be verified.

It is of course too early to assess whether the new algorithms will turn out
to be of practical interest. The appraisal of their numerical behaviour
is the object of ongoing research.

{\footnotesize
  
%\bibliography{/home/pht/bibs/refs}
\bibliographystyle{plain}

}

\appendix

\appnumsection{Appendix: technical lemmas}

\llem{wlogw}{Let $w >0$ and suppose that
  \beqn{w-implicit}
  w^\alpha \leq \beta \log(2w).
  \eeqn
  for some $\alpha\in (0,1)$ and $\beta$ such that
  \beqn{beta-cond}
  \beta > \frac{3\alpha}{2^\alpha}.
  \eeqn
  Then
  \beqn{w-explicit}
  w \leq \sigma(\alpha,\beta)
  \eqdef
  \left[-\frac{\beta}{\alpha}W_{-1}\left(-\frac{\alpha}{\beta\,2^\alpha}\right) \right]^\sfrac{1}{\alpha}
  \eeqn
  where $W_{-1}(\cdot)$ is the second branch of the Lambert function \cite{Corletal96}.
}

\proof{
  First note that \req{w-implicit} is equivalent to
  \[
  \frac{1}{2^\alpha}(2 w)^\alpha \leq
  \frac{\beta}{\alpha}\log\Big((2 w)^\alpha \Big)
  \]
  Setting now $u = (2 w)^\alpha$, one obtains that
  \beqn{vt-ineq}
  \omega(u)\eqdef \frac{1}{2^\alpha}u - \frac{\beta}{\alpha} \log(u)\leq 0.
  \eeqn
  But $\omega(u)$ is convex for $u>0$ and tends to infinity if $u$ tends to
  zero or to infinity.  Moreover, it achieves its minimum at $u_{\min}=
  \beta 2^\alpha/\alpha$, at which it takes the value
  \[
  \omega(u_{\min}) =
  \frac{\beta}{\alpha}\left(1-\log\left(\frac{\beta \, 2^\alpha}{\alpha}\right)\right)
  < 0,
  \]
  where the inequality results from \req{beta-cond}.  Hence $\omega(u)$ has
  two real roots $u_1\leq u_2$ and the set of $u$ for
  which \req{vt-ineq} holds is bounded above by $u_2$.  By definition,
  \[
  \log(u_2)-\frac{\alpha}{\beta\, 2^\alpha}\, u_2 = 0,
  \]
  which is
  \[
  u_2 e^{-\frac{\alpha}{\beta\, 2^\alpha}\,u_2} = 1.
  \]
  Defining now $z = -\frac{\alpha}{\beta\,2^\alpha}\,u_2$, we obtain
  that
  \[
  z e^z =  -\frac{\alpha}{\beta\, 2^\alpha} .
  \]
  By definition of the Lambert function, this gives that
  \[
  u_2 = -\frac{\beta \, 2^\alpha}{\alpha}z
  = -\frac{\beta 2^\alpha}{\alpha}
  W_{-1}\left(-\frac{\alpha}{\beta \, 2^\alpha}\right) > 0
  \]
  which is well-defined because \req{beta-cond} implies that
  $-\frac{\alpha}{\beta\,2^\alpha}\in [-\frac{1}{e},0)$.
  Since $w = u^\sfrac{1}{\alpha}/2$, this implies \req{w-explicit}.
} %epr

\llem{techa}{Let $a\geq 0$ and $b\geq0$.
Suppose that, for some $\mu \in (0,\half)$, $\nu\in(0,\third)$ and
some $\theta_{a,0}, \theta_{b,0}$ and $\theta_0 \geq 0$,
\beqn{techa:hyp}
\theta_{a,0}a^{1-\mu} + \theta_{b,0}b^{1-2\nu} \leq  \theta_0.
\eeqn
Then
\[
a \leq \left(\frac{\theta_0}{\theta_{a,0}}\right)^\sfrac{1}{1-\mu}
\tim{ and }
b \leq \left(\frac{\theta_0}{\theta_{b,0}}\right)^\sfrac{1}{1-2\nu}.
\]
}

\proof{Obvious from the inequalities $\theta_{a,0}a^{1-\mu} \leq \theta_{a,0}a^{1-\mu}+ \theta_{b,0}b^{1-2\nu}$
and $\theta_{b,0}b^{1-2\nu}\leq \theta_{a,0}a^{1-\mu} + \theta_{b,0}b^{1-2\nu}$.
} % epr

\llem{techb}{
Let $a\geq 0 $ and $b\geq 0$. Suppose that, for some $\mu \in (0,\half)$,
some $\theta_{a,1} > 0$ and some $\theta_1 \geq 0$,
\beqn{techb:hyp-a}
a^{1-\mu} \leq \theta_{a,1} a^{1-2\mu} + \theta_1.
\eeqn
Then
\[
a \leq \max\left[ (2\theta_1)^\sfrac{1}{1-\mu},(2\theta_{a,1})^\sfrac{1}{\mu}\right].
\]
Symmetrically, if $\nu\in(0,\third)$, $\theta_{b,1} > 0$ and
\beqn{techb:hyp-b}
b^{1-2\nu} \leq \theta_{b,1} b^{1-3\nu} + \theta_1.
\eeqn
Then
\[
b \leq \max\left[ (\theta_1)^\sfrac{1}{1-2\nu},(2\theta_{b,1})^\sfrac{1}{\nu}\right].
\]
}

\proof{
  Suppose first that $\theta_{a,1} a^{1-2\mu}\leq \theta_1$. Then
  $a^{1-\mu}\leq 2\theta_1$ and thus $a\leq (2\theta_1)^\sfrac{1}{1-\mu}$
 Suppose now that $\theta_{a,1} a^{1-2\mu}> \theta_1$.
Then $a^{1-\mu}\leq 2\theta_{a,1} a^{1-2\mu}$, that is $a\leq (2\theta_a)^\sfrac{1}{\mu}$.
The proof of the second part is similar.
} %epr

\llem{techc}{
Let $a\geq 0$ and $b\geq0$. Suppose that, for some $\mu \in (0,\half)$, $\nu\in(0,\third)$ and
some $\theta_a, \theta_b > 0$ and $\theta_2 \geq 0$,
\beqn{techc:hyp}
a^{1-\mu} + b^{1-2\nu} \leq \theta_{a,2} a^{1-2\mu} + \theta_{b,2} b^{1-3\nu} + \theta_2.
\eeqn
Then
\[
a \leq \max\left[
  \left(\frac{\theta_2}{\theta_{a,2}}\right)^\sfrac{1}{1-2\mu},
  2^\frac{1}{1-\mu}\left(2\theta_{b,2}\right)^\sfrac{1-2\nu}{\nu(1-\mu)},
  \left(4\theta_{a,2}\right)^\sfrac{1}{\mu}
  \right]
\]
and
\[
b \leq
\max\left[
  \left(\frac{\theta_2}{\theta_{b,2}}\right)^\sfrac{1}{1-3\nu},
  2^\frac{1}{1-2\nu}\left(2\theta_{a,2}\right)^\sfrac{1-\mu}{\mu(1-2\nu)},
  \left(4\theta_{b,2}\right)^\sfrac{1}{\nu}
  \right].
\]
}

\proof{%%%
  Suppose first that
  \beqn{t1:c1}
  \theta_{a,2} a^{1-2\mu}+\theta_{b,2} b^{1-3\nu} \leq \theta_2.
  \eeqn
  Then, from Lemma~\ref{techa},
  \beqn{t1:ab1}
  a \leq \left(\frac{\theta_2}{\theta_{a,2}}\right)^\sfrac{1}{1-2\mu}
  \tim{ and }
  b \leq \left(\frac{\theta_2}{\theta_{b,2}}\right)^\sfrac{1}{1-3\nu}.
  \eeqn
  Suppose now that \req{t1:c1} fails, and thus that
  \beqn{t1:e1}
  \theta_{a,2} a^{1-2\mu}+\theta_{b,2} b^{1-3\nu} + \theta_2
  \leq 2\theta_{a,2} a^{1-2\mu}+2\theta_{b,2} b^{1-3\nu}.
  \eeqn
  Assume also that
  \beqn{abs1}
  a > \left(2\theta_{a,2}\right)^{\sfrac{1}{\mu}}
  \tim{and}
  b > \left(2\theta_{b,2}\right)^{\sfrac{1}{\nu}}
  \eeqn
  Then,
  \[
  2\theta_{a,2} a^{1-2\mu} + 2\theta_{b,2} b^{1-3\nu}
  <  a^{1-\mu} +  b^{1-2\nu}
  \]
  and so, using \req{techc:hyp} and \req{t1:e1},
  \[
  a^{1-\mu} + b^{1-2\nu}
  \leq \theta_{a,2}a^{1-2\mu}+\theta_{b,2} b^{1-3\nu}+\theta_2
  < a^{1-\mu} +  b^{1-2\nu},
  \]
  which is impossible.  Hence \req{abs1} cannot hold, and at least one of its
  inequalities must fail.  Suppose that it is the first, that is
  \beqn{ee1}
  a \leq \left(2\theta_{a,2}\right)^{\sfrac{1}{\mu}} \eqdef \kappa_1.
  \eeqn
  Then \req{techc:hyp} and \req{t1:e1} give that
  \[
  b^{1-2\nu}
  \leq a^{1-\mu}+b^{1-2\nu}
  \leq 2\theta_{a,2}\kappa_1^{1-2\mu} + 2\theta_{b,2} b^{1-3\nu}
  \]
  and we may apply Lemma~\ref{techb} with $\theta_{b,1} = 2\theta_{b,2}$ and
  $\theta_1 = 2\theta_{a,2}\kappa_1^{1-2\mu}$ to deduce that
  \[
  b \leq
  \max\left[\left(4\theta_{a,2}\kappa_1^{1-2\mu}\right)^{\sfrac{1}{1-2\nu}},
            \left(4\theta_{b,2}\right)^{\sfrac{1}{\nu}}\right].
  \]
  Symmetrically, we deduce that if the second inequality of \req{abs1} fails,
  that is if
  \[
  b \leq \left(2\theta_{b,2}\right)^{\sfrac{1}{\nu}} \eqdef \kappa_2,
  \]
  then, applying Lemma~\ref{techb} with $\theta_{a,1}= 2 \theta_{a,2}$
  and $\theta_1 = 2\theta_{b,2}\kappa_2^{1-3\nu}$,
  \[
  a \leq
  \max\left[\left(4\theta_{b,2}\kappa_2^{1-3\nu}\right)^{\sfrac{1}{1-\mu}},
            \left(4\theta_{a,2}\right)^{\sfrac{1}{\mu}}\right].
  \]
  Combining the two cases yields the desired result.
} %epr

\llem{techd}{
Let $a\geq 0$ and $b \geq 0$. Suppose that, for some $\mu \in (0,\half)$, $\nu\in(\third,1)$ and
some $\theta_{a,3}>0$, $\theta_3 \geq 0$,
\beqn{techd:hyp}
a^{1-\mu} + b^{1-2\nu} \leq \theta_{a,3} a^{1-2\mu}+\theta_3.
\eeqn
Then
\[
a \leq \max\left[
  (2\theta_3)^\sfrac{1}{1-\mu},
  (2\theta_{a,3})^\sfrac{1}{\mu}
  \right] = \kappa_{a,3}
\tim{ and }
b \leq \left(\theta_{a,3} \kappa_{a,3}^{1-2\mu}+\theta_3\right)^\sfrac{1}{1-2\nu}.
\]
Symmetrically, if $\theta_{b,3}>0$ and
\[
a^{1-\mu} + b^{1-2\nu} \leq \theta_{b,3} b^{1-3\nu}+\theta_3,
\]
then
\[
b \leq \max\left[
  (2\theta_3)^\sfrac{1}{1-2\nu},
  (2\theta_{b,3})^\sfrac{1}{\nu}
  \right] = \kappa_{b,3}
\tim{ and }
a \leq \left(\theta_{b,3} \kappa_{b,3}^{1-3\nu}+\theta_3\right)^\sfrac{1}{1-\mu}.
\]
}

\proof{
From \req{techd:hyp}, we have that
\[
a^{1-\mu}
\leq a^{1-\mu} + b^{1-2\nu}
\leq \theta_{a,3} a^{1-2\mu}+\theta_2
\]
and we may apply Lemma~\ref{techb} with $\theta_{a,1}=\theta_{a,3}$ and $\theta_1=\theta_3$
to deduce that
\[
a \leq \max\left[ \left(2\theta_3\right)^\sfrac{1}{1-\mu}, (2\theta_{a,3})^\sfrac{1}{\mu}\right]
\eqdef \kappa_a
\]
From the inequality $b^{1-2\nu} \leq a^{1-\mu}+b^{1-2\nu}$ and
\req{techd:hyp}, we also obtain that 
\[
b \leq \left(\theta_{a,3}
\kappa_a^{1-2\mu}+\theta_3\right)^\sfrac{1}{1-2\nu}.
\]
} % epr

\llem{teche}{
Let $a>0$ and $b>0$. Suppose that, for some $\nu\in(0,\third]$,
some $\theta_{a,4} \geq 1$ and some $\theta_4\geq0$,
\beqn{teche:hyp}
a^\half + b^{1-2\nu} \leq \theta_{a,4} \log(2a) +\theta_4.
\eeqn
Then
\[
a \leq
\max\left[\frac{1}{2}e^\sfrac{\theta_4}{\theta_{a,4}}, \sigma\big(\half,2\theta_{a,4}\big)\right]
= \kappa_{a,4}
\tim{ and }
b \leq \Big(\theta_{a,4} \log(2\kappa_{a,4})+\theta_4\Big)^\sfrac{1}{1-2\nu}.
\]
Symmetrically, if  $\theta_{b,4}\geq 1$, $\mu\in(0,\half]$ and
\[
a^{1-\mu} + b^\third \leq \theta_{b,4} \log(2b) +\theta_4,
\]
then
\[
b \leq
\max\left[\frac{1}{2}e^\sfrac{\theta_4}{\theta_{b,4}}, \sigma\big(\half,2\theta_{b,4}\big)\right]
= \kappa_{b,4}
\tim{ and }
a \leq \Big(\theta_{b,4} \log(2\kappa_{b,4})+\theta_4)\Big)^\sfrac{1}{1-\mu}.
\]
}

\proof{
Suppose first that $\theta_{a,4} \log(2a)\leq \theta_4$.  Then
\beqn{te:b1}
a \leq \frac{1}{2}e^\sfrac{\theta_4}{\theta_{a,4}}.
\eeqn
Otherwise, \req{teche:hyp} gives that
\[
a^\half \leq a^\half + b^{1-2\nu} \leq 2\theta_{a,4} \log(2a)
\]
from which one deduces using Lemma~\ref{wlogw} with $\alpha = \half$ and $\beta=2\theta_{a,4}$
(which is allowed since $2\theta_{a,4}\geq 2> 3/2^\sfrac{5}{2}$ implies \req{beta-cond}) that
\[
a \leq \sigma(\half, 2\theta_{a,4}),
\]
where $\sigma(\cdot,\cdot)$ is defined in \req{w-explicit}.
This inequality and \req{te:b1} give the desired bound on $a$. Substituting
this in \req{teche:hyp} gives the bound on $b$. 
The proof of the symmetric statement is similar, in which the use of
Lemma~\ref{wlogw} is now allowed because $\theta_{b,4} \geq 1 > 1/2^\sfrac{4}{3}$
again implies \req{beta-cond}.
} %epr

\llem{techf}{
Let $a>0$ and $b\geq0$. Suppose that, for some $\nu\in(0,\third)$,
some $\theta_{a,5} \geq 1$, $\theta_{b,5}>0$ and some $\theta_5\geq0$,
\beqn{techf:hyp}
a^\half + b^{1-2\nu} \leq \theta_{a,5} \log(2a) + \theta_{b,5} b^{1-3\nu} +\theta_5,
\eeqn
Then
\[
a \leq \max\left[
  \frac{1}{2}e^\sfrac{\theta_5}{\theta_{a,5}},
  \sigma(\half,4\theta_{a,5})
  \frac{1}{2}e^\sfrac{{\theta_{b,5}(2\theta_{b,5})^\sfrac{1-3\nu}{\nu}}}{\theta_{a,5}}
  \right]
\]
and
\[
b \leq \max\left[
  \left(\frac{\theta_5}{\theta_{b,5}}\right)^\sfrac{1}{1-3\nu},
  \left(2\theta_{a,5}\log(2\sigma_a)\right)^\sfrac{1}{1-2\nu},
  \Big(4\theta_{b,5}\Big)^\frac{1}{\nu},
   (2\theta_{b,5})^\frac{1}{\nu}
  \right].
\]
with $\sigma_a= \sigma\big(\half,2\theta_{a,5}\big)$.
Symmetrically, if $\theta_{a,5} \geq 1$, $\mu\in (0,\half)$, $b>0$  and
\[
a^{1-\mu} + b^\third
\leq \theta_{a,5} a^{1-2\mu} + \theta_{b,5} \log(2b) +\theta_5,
\]
then
\[
a \leq \max\left[
  \left(\frac{\theta_5}{\theta_{a,5}}\right)^\sfrac{1}{1-2\mu},
  \left(2\theta_{b,5}\log(2\sigma_b)\right)^\sfrac{1}{1-\mu},
  \Big(4\theta_{a,5}\Big)^\frac{1}{\mu},
   (2\theta_{a,5})^\frac{1}{\mu}
  \right]
\]
and
\[
b \leq \max\left[
  \frac{1}{2}e^\sfrac{\theta_5}{\theta_{b,5}},
  \sigma(\half,4\theta_{b,5})
  \frac{1}{2}e^\sfrac{{\theta_{a,5}(2\theta_{a,5})^\sfrac{1-2\mu}{\mu}}}{\theta_{b,5}}
  \right]
\]
with $\sigma_b= \sigma\big(\half,2\theta_{b,5}\big)$.
}

\proof{
Suppose first that,
\beqn{tf:c1}
\theta_{a,5}\log(2a) + \theta_{b,5}b^{1-3\nu} \leq \theta_5.
\eeqn
Then
\beqn{tf:b1}
a \leq \frac{1}{2}e^{\sfrac{\theta_5}{\theta_{a,5}}}
\tim{ and }
b\leq \left(\frac{\theta_5}{\theta_{b,5}}\right)^\sfrac{1}{1-3\nu}.
\eeqn
Suppose now that \req{tf:c1} fails. Then, from \req{techf:hyp},
\beqn{tf:c2}
a^\half + b^{1-2\nu} \leq 2\theta_{a,5} \log(2a) + 2\theta_{b,5} b^{1-3\nu}.
\eeqn
If 
\beqn{abs5}
a > \sigma\big(\half,2\theta_{a,5}\big)
\tim{ and }
b > (2\theta_{b,5})^\sfrac{1}{\nu},
\eeqn
we obtain, using Lemma~\ref{wlogw} (which we may apply because $\theta_{a,5}
\geq 1 > 3/2^\sfrac{5}{2}$), \req{tf:c2} and \req{techf:hyp} that
\[
a^\half + b^{1-2\nu}
\leq 2\theta_{a,5} \log(2a) + 2\theta_{b,5} b^{1-3\nu}
<a^\half+b^{1-2\nu},
\]
which is impossible. Hence one of the inequalities of \req{abs5} must be
violated.  Suppose that
\[
a\leq \sigma\big(\half,2\theta_{a,5}\big) \eqdef \sigma_a.
\]
Using Lemma~\ref{wlogw} again and \req{tf:c2}, this implies that
\[
b^{1-2\nu}
\leq a^\half+b^{1-2\nu}
\leq 2\theta_{a,5}\log(2\sigma_a) + 2\theta_{b,5} b^{1-3\nu}
\]
and we deduce from Lemma~\ref{techb} with $\theta_{b,1}=2\theta_{b,5}$
and $\theta_1= 2\theta_{a,5}\log(2\sigma_a)$ that
\[
b \leq
\max\left[\left(2\theta_{a,5}\log(2\sigma_a)\right)^\sfrac{1}{1-2\nu},
  \Big(4\theta_{b,5}\Big)^\frac{1}{\nu}\right].
\]
If we now suppose that $b \leq (2\theta_{b,5})^\sfrac{1}{\nu}$, then
\req{tf:c2} ensures that
\[
a^\half
\leq a^\half+b^{1-2\nu}
\leq 2\theta_{a,5} \log(2a)  + 2\theta_{b,5} (2\theta_{b,5})^\sfrac{1-3\nu}{\nu},
\]
and we now obtain from Lemma~\ref{teche} with
$\theta_{a,4} = 2\theta_{a,5}$ and $\theta_4=2\theta_{b,5} (2\theta_{b,5})^\sfrac{1-3\nu}{\nu}$
that
\[
a \leq
\max\left[
  \frac{1}{2}e^\sfrac{\theta_{b,5}(2\theta_{b,5})^\sfrac{1-3\nu}{\nu}}{\theta_{a,5}},
  \sigma\big(\half,4\theta_{a,5}\big)
  \right].
\]
} %epr

\llem{techg}{
Let $a>0$ and $b>0$. Suppose that, for 
some $\theta_{a,6} \geq 1$, $\theta_{b,6} \geq 1$ and some $\theta_6\geq0$,
\beqn{techg:hyp}
a^\half + b^\third
\leq \theta_{a,6} \log(2a) + \theta_{b,6} \log(2b) +\theta_6,
\eeqn
where $2a\geq\varsigma$ and $2b\geq\varsigma$. Then
\[
a \leq \max\left[
  \frac{1}{2}e^\sfrac{\theta_6 + |\log(\varsigma)|}{\theta_{a,6}},
  \sigma(\half,2\theta_{a,6}),
  \sigma(\third,2\theta_{b,6})e^\sfrac{\theta_{b,6}}{2\theta_{a,6}},
  \sigma(\third,4\theta_{a,6})
  \right]
\]
and
\[
b \leq \max\left[
  \frac{1}{2}e^\sfrac{\theta_6 + |\log(\varsigma)|}{\theta_{b,6}},
  \sigma(\half,2\theta_{b,6}),
  \sigma(\third,2\theta_{ba6})e^\sfrac{\theta_{a,6}}{2\theta_{b,6}},
  \sigma(\third,4\theta_{b,6})
  \right].
\]
}

\proof{%%%%
  Suppose first that
  \beqn{tg:c1}
  \theta_{a,6} \log(2a) + \theta_{b,6} \log(2b) \leq \theta_6.
  \eeqn
  Then 
  \[
  \theta_{a,6} \log(2a)  \leq \theta_6 + |\log(\varsigma)|
  \tim{ and }
  \theta_{b,6} \log(2b)  \leq \theta_6 + |\log(\varsigma)|
  \]
  and hence
  \beqn{tg:b1}
  a \leq \frac{1}{2}e^\sfrac{\theta_6 + |\log(\varsigma)|}{\theta_{a,6}}
  \tim{ and }
  b \leq \frac{1}{2}e^\sfrac{\theta_6 + |\log(\varsigma)|}{\theta_{b,6}}.
  \eeqn
  Suppose now that \req{tg:c1} fails, and thus \req{techg:hyp} implies that
  \beqn{tg:c2}
  a^\half + b^\third
  \leq 2 \theta_{a,6} \log(2a) + 2 \theta_{b,6} \log(2b).
  \eeqn
  Assume also that
  \beqn{abs6}
  a > \sigma\big(\half,2\theta_{a,6}\big)
  \tim{ and }
  b > \sigma\big(\third,2\theta_{b,6}\big).
  \eeqn
  Then, using \req{techg:hyp} and \req{tg:c2},
  \[
  a^\half + b^\third
  \leq 2 \theta_{a,6} \log(2a) + 2 \theta_{b,6} \log(2b)
  < a^\half + b^\third,
  \]
  which is impossible.  Hence one of the inequalities of \req{abs6} must fail.
  If
  \[
  a \leq \sigma\big(\half,2\theta_{a,6}\big),
  \]
  then \req{techg:hyp} gives that
  \[
  b^\third \leq a^\half+ b^\third
  \leq \theta_{a,6}\log\Big(2\sigma\big(\half,2\theta_{a,6}\big)\Big)+ 2\theta_{b,6} \log(2b).
  \]
  and Lemma~\ref{teche} with $\theta_{b,4}= 2\theta_{b,6}$ and
  $\theta_4 = \theta_{a,6}\log\Big(2\sigma\big(\half,2\theta_{a,6}\big)\Big)$ then implies
  that
  \[
  b\leq
  \max\left[\frac{1}{2}e^\sfrac{\theta_{a,6}\log(2\sigma(\half,2\theta_{a,6}))}{2\theta_{b,6}},
    \sigma\big(\half,4\theta_{b,6}\big)\right]
  =\max\left[\sigma\big(\half,2\theta_{a,6}\big)e^\sfrac{\theta_{a,6}}{2\theta_{b,6}},
    \sigma\big(\half,4\theta_{b,6}\big)\right].
  \]
  Symmetrically, if
  \[
  b < \sigma\big(\third,2\theta_{b,6}\big),
  \]
  then
  \[
  a\leq
  \max\left[\sigma(\third,2\theta_{b,6})e^\sfrac{\theta_{b,6}}{2\theta_{a,6}},
    \sigma\big(\third,4\theta_{a,6}\big)\right].
  \]
}% epr

\end{document}